\documentclass[a4paper, 11pt, twoside]{article}
\usepackage{CJK}
\usepackage{latexsym,bm}
\usepackage[tbtags]{amsmath}
\usepackage{algorithm}
\usepackage{algorithmic}
\usepackage{amssymb}
\usepackage{fancyhdr}
\usepackage{cite}
\usepackage{tabularx}
\usepackage{booktabs}
\usepackage{dcolumn}
\usepackage{graphicx}
\usepackage{wrapfig}
\usepackage{multicol}
\usepackage{verbatim}
\usepackage{mathrsfs}
\usepackage{multirow}

\usepackage{graphicx}

\def\nkyear{}                           
\def\nkvol{}                              
\def\nkno{}                                
\def\firstpage{1}                           
\def\DOI{}
\def\shorttitle{A new stable and avoiding inversion iteration for computing matrix square root}
\def\shortauthor{Li ZHU, Keqi YE, Yuelin ZHAO, Feng WU, Jiqiang HU and Wanxie ZHONG }
\usepackage{Shangda1226}                    

\newcommand{\supercite}[1]{\!\!\textsuperscript{\cite{#1}}} 
\def\ge{\geqslant}\def\le{\leqslant}

\begin{document} 

\begin{CJK*}{GBK}{song}
\title{{\large  \textbf{A new stable and avoiding inversion iteration for computing matrix square root}}
\thanks{\dag\,\,Corresponding author, E-mail: vonwu@dlut.edu.cn\;
\newline        
Project supported by the National Natural Science Foundation of
China grants (Nos. 11472076 and 51609034), the Science Foundation of Liaoning Province of China (No. 2021-MS-119), the Dalian Youth Science and Technology Star project (No. 2018RQ06),
and the Fundamental Research Funds for the Central Universities grant (Nos. DUT20RC(5)009 and DUT20GJ216).}}
\author{\small{Li ZHU$^1$,\quad Keqi YE$^1$,\quad Yuelin ZHAO$^1$,} 
\\
\small{Feng WU$^{1,\dag}$,\quad Jiqiang HU$^1$,\quad Wanxie ZHONG$^1$}
\\[2mm]
\footnotesize{1. Key Laboratory of Structural Analysis of Industrial Equipment, Department of Mechanics,}
\\
\footnotesize{Faculty of Vehicle Engineering and Mechanics, Dalian University of Technology,}
\\
\footnotesize{Dalian 116024, People's Republic of China;}
\\[2mm]}

\maketitle 
\thispagestyle{first} \footnotesize
\begin{abstract}
\noindent \textbf{Abstract~~~}The objective of this research was to compute the principal matrix square root with sparse approximation. A new stable iterative scheme avoiding fully matrix inversion (SIAI) is provided. The analysis on the sparsity and error of the matrices involved during the iterative process is given. Based on the bandwidth and error analysis, a more efficient algorithm combining the SIAI with the filtering technique is proposed. The high computational efficiency and accuracy of the proposed method are demonstrated by computing the principal square roots of different matrices to reveal its applicability over the existing methods.
\\[2mm]
\textbf{Key words~~~}Matrix square root, Iterative algorithm, Error analysis, Bandwidth analysis
\\[2mm]
\textbf{Chinese Library Classification~~~}O29
\\
\textbf{2020 Mathematics Subject Classification~~~}65F45
\end{abstract}

\section{Introduction}\label{S1}
\small Matrix square root is one of the most commonly involved matrix functions. Unlike solving the square root of a scalar, a matrix may not have a square root. And sometimes it is possible for a matrix to have two or more square roots. For example, the non-singular Jordan block has two square roots, and any involutory matrix is one square root of the identity matrix. If the matrix $\mathbf{A}$ is singular, then the existence of its square root depends on the structure of Jordan blocks corresponding to the zero eigenvalues\supercite{1}. If the matrix $\mathbf{A}$ is a real matrix, it probably does not have a real square root. The theory on the matrix square root is very complex, and different theories have been developed for different matrices with different properties. 

It has been demonstrated that for any matrix $\mathbf{A}\in {{\mathbb{C}}^{n\times n}}$ with no nonpositive real eigenvalues, there is only one square root $\mathbf{X}$, satisfying that the real parts of its eigenvalues are all positive. Such matrix $\mathbf{X}$ is denoted by ${{\mathbf{A}}^{1/2}}$, and is called the principal square root\supercite{2}. The  principal matrix square root plays an important role in many math problems, such as definite generalized eigenvalue problem\supercite{3}, polar decomposition\supercite{4}, geometric mean value\supercite{5}, matrix sign function\supercite{6} and matrix logarithm function\supercite{7}. The matrix square root also plays an important role in many practical problems, such as dynamics problem, deep learning\supercite{8} and machine vision$^{[9\textrm{--}11]}$.

Many algorithms have been developed for solving the principal square root of a matrix. One of the most famous algorithms is the Newton iteration\supercite{12}
\begin{equation}\label{E2}
\left\{ \begin{aligned}
  & {{\mathbf{X}}_{k}}{{\mathbf{H}}_{k}}+{{\mathbf{H}}_{k}}{{\mathbf{X}}_{k}}=\mathbf{A}-\mathbf{X}_{k}^{2}, \\ 
 & {{\mathbf{X}}_{k+1}}={{\mathbf{X}}_{k}}+{{\mathbf{H}}_{k}},\\ 
\end{aligned} \right.\ k=0,1,\cdots.
\end{equation}
Choosing a suitable initial value ${{\mathbf{X}}_{0}}$, a sequence $\left\{ {{\mathbf{X}}_{k}} \right\}$ that converges to ${{\mathbf{A}}^{1/2}}$ can be given. This method is numerically stable\supercite{8}, but it requires solving a Sylvester equation in each iteration, which is expensive. If we set ${{\mathbf{X}}_{0}}=\mathbf{A}$, then ${{\mathbf{X}}_{k}}$ commutes with ${{\mathbf{H}}_{k}}$, and (\ref{E2}) can be rewritten as
\begin{equation}\label{E3}
{{\mathbf{X}}_{k+1}}=\frac{1}{2}\left( {{\mathbf{X}}_{k}}+\mathbf{AX}_{k}^{-1} \right),\ {{\mathbf{X}}_{0}}=\mathbf{A}.
\end{equation}
However, Higham\supercite{12} pointed out that the iterative scheme (\ref{E3}), despite avoiding solving the Sylvester equations in (\ref{E2}), is unstable in most cases.

To solve the stability problem of (\ref{E3}), some modified versions have been proposed, such as the Denman and Beavers (DB) method$^{[13\textrm{--}14]}$, and the Meini iteration based on the cyclic reduction algorithm (CR)\supercite{15}. The DB method is based on the iteration of matrix sign function, and can be expressed as:
\begin{equation}\label{E4}
\left\{ \begin{aligned}
  & {{\mathbf{X}}_{0}}=\mathbf{A},{{\mathbf{Y}}_{0}}=\mathbf{I}, \\ 
 & {{\mathbf{X}}_{k+1}}=\frac{1}{2}\left( {{\mathbf{X}}_{k}}+\mathbf{Y}_{k}^{-1} \right), \\ 
 & {{\mathbf{Y}}_{k+1}}=\frac{1}{2}\left( {{\mathbf{Y}}_{k}}+\mathbf{X}_{k}^{-1} \right), \\ 
\end{aligned} \right.\ k=0,1,\cdots,
\end{equation}
which generates a sequence $\left\{ {{\mathbf{X}}_{k}} \right\}$ that converges to ${{\mathbf{A}}^{1/2}}$ and a sequence $\left\{ {{\mathbf{Y}}_{k}} \right\}$ that converges to ${{\mathbf{A}}^{-1/2}}$. The scheme of the CR method can be defined by
\begin{equation}\label{E5}
\left\{ \begin{aligned}
  & {{\mathbf{Z}}_{0}}=2\left( \mathbf{I}+\mathbf{A} \right),{{\mathbf{Y}}_{0}}=\mathbf{I}-\mathbf{A}, \\ 
 & {{\mathbf{Y}}_{k+1}}=-{{\mathbf{Y}}_{k}}\mathbf{Z}_{k}^{-1}{{\mathbf{Y}}_{k}}, \\ 
 & {{\mathbf{Z}}_{k+1}}={{\mathbf{Z}}_{k}}-2{{\mathbf{Y}}_{k}}\mathbf{Z}_{k}^{-1}{{\mathbf{Y}}_{k}}, \\ 
\end{aligned} \right.k=0,1,\cdots.
\end{equation}
The sequence $\left\{ {{\mathbf{Z}}_{k}} \right\}$ generated by (\ref{E5}) converges to $\frac{1}{4}{{\mathbf{A}}^{1/2}}$. In Ref.\supercite{16}, Iannazzo derived an iteration based on Newton's method,
\begin{equation}\label{E6}
\left\{ \begin{aligned}
  & {{\mathbf{X}}_{0}}=\mathbf{A},{{\mathbf{H}}_{0}}=\frac{1}{2}\left( \mathbf{I}-\mathbf{A} \right), \\ 
 & {{\mathbf{X}}_{k+1}}={{\mathbf{X}}_{k}}+{{\mathbf{H}}_{k}}, \\ 
 & {{\mathbf{H}}_{k+1}}=-\frac{1}{2}{{\mathbf{H}}_{k}}\mathbf{X}_{k+1}^{-1}{{\mathbf{H}}_{k}}. \\ 
\end{aligned} \right.\ k=0,1,\cdots.
\end{equation}
which is equivalent to (\ref{E5}) but generates a sequence $\left\{ {{\mathbf{X}}_{k}} \right\}$ that converges directly to ${{\mathbf{A}}^{1/2}}$.

The above algorithms (DB, CR, IN) have been proved to be stable$^{[15\textrm{--}17]}$. However, in the face of large-scale sparse matrices, these methods are time-consuming because they all require matrix inversions more or less, which cost a lot of computational effort. How to design a stable scheme avoiding matrix inversion is obviously a valuable issue, and will be discussed in this paper.

	Another limitation on computing the root of a large sparse matrix is the memory required, as the root of a large sparse matrix is often dense or full. Demko et al.\supercite{18}, while studying the inverse of symmetric positive definite banded matrices, found that the absolute values of the elements of the final inverse matrix decay exponentially as they move away from the main diagonal, and the more diagonally dominant the matrix is, the stronger the decay. Benzi et al.\supercite{19} successfully extended this property to general matrix functions and proved that the elements of the matrix functions of all symmetric banded matrices move away from the main diagonal in an exponentially decaying manner. In Ref.\supercite{20}, they also summarized this decay property as a localization phenomenon and proposed that a part of the minimal elements can be neglected in numerical calculations without affecting the computational error, which has a positive impact on the design of fast approximation algorithms. Inspired by this, Gao et al.\supercite{21} proposed a filtering technique and used it for the solution of large-scale sparse matrix exponentials, showing great efficiency. Wu et al.\supercite{22} introduced the $\varepsilon$-bandwidth to measure the sparsity of a matrix. If the $\varepsilon$-bandwidth of a matrix is far smaller than its dimension, the matrix is treated as a nearly sparse matrix. For the matrix whose principal square root matrix is nearly sparse matrix, it may be a potential way to combine the filtering technique with the existing iterative scheme to construct an efficient method for the computation of the square root of large and sparse matrices.

This paper is concerned with the numerical computation of the principal square root of a large-scale sparse matrix. In Section 2, a new stable iteration avoiding inversions (SIAI) is developed and a stability proof is presented. In Section 3, we pointed out the applicability of the filtering technique by using the bandwidth analysis. In Section 4, the filtering algorithm is introduced into the proposed SIAI, and the adaptive selection of the filtering threshold is discussed in detail. Numerical experiments are given in Section 5 to test the accuracy and efficiency of the proposed method. Some conclusions are summarized in Section 6. 

\section{A stable iteration avoiding inversions}\label{S2}

To study the principal square root of matrix $\mathbf{A}\in {{\mathbb{C}}^{n\times n}}$ is to solve the following equation
\begin{equation}\label{E7}
{{\mathbf{X}}^{2}}-\mathbf{A}=\mathbf{0}.
\end{equation}
Applying Newton' s method to this equation will introduce matrix inversions more or less. For the large-scale sparse matrix, the inverse matrix is usually dense which makes the computation very expensive. To avoid or reduce matrix inversions, Ref.\supercite{23} recommended to consider the following equation
\begin{equation}\label{E8}
{{\mathbf{X}}^{-2}}-{{\mathbf{A}}^{-1}}=\mathbf{0}.
\end{equation}
Applying Newton's method to (\ref{E8}) will yield the following iteration format:
\begin{equation}\label{E9}
\left\{ \begin{aligned}
  & {{\mathbf{X}}_{k+1}}={{\mathbf{X}}_{k}}\left( \mathbf{I}+\frac{1}{2}{{\mathbf{Y}}_{k}} \right), \\ 
 & {{\mathbf{Y}}_{k+1}}=\mathbf{I}-{{\mathbf{A}}^{-1}}\mathbf{X}_{k+1}^{2},\ k\ge 0. \\ 
\end{aligned} \right.
\end{equation}
Actually, the above iteration is a variant of the inverse Newton iteration\supercite{24}. If replace ${{\mathbf{A}}^{-1}}$ with $\mathbf{A}$, there will be ${{\mathbf{X}}_{k}}\to {{\mathbf{A}}^{-0.5}}$ as $k\to \infty $.

\subsection{How to choose a suitable ${{\mathbf{X}}_{0}}$}\label{S2.1}

According to the Ref.\supercite{25}, if $\mathbf{A}{{\mathbf{X}}_{0}}={{\mathbf{X}}_{0}}\mathbf{A}$, one has
\begin{equation}\label{E10}
{{\mathbf{Y}}_{k+1}}=\frac{3}{4}\mathbf{Y}_{k}^{2}+\frac{1}{4}\mathbf{Y}_{k}^{3}.
\end{equation}
Then, taking the norm of both sides of (\ref{E10}) yields
\begin{equation}\label{E11}
\left\| {{\mathbf{Y}}_{k+1}} \right\|\le \frac{3}{4}{{\left\| {{\mathbf{Y}}_{k}} \right\|}^{2}}+\frac{1}{4}{{\left\| {{\mathbf{Y}}_{k}} \right\|}^{3}}
\end{equation}
for any consistent norm. Therefore, if we can find such a ${{\mathbf{X}}_{0}}$ which commutes with $\mathbf{A}$ and satisfies $\left\| {{\mathbf{Y}}_{0}} \right\|<1$, then the iteration (\ref{E9}) will generate the sequence $\left\{ {{\mathbf{X}}_{k}} \right\}$ converging to ${{\mathbf{A}}^{1/2}}$. The most common selection is ${{\mathbf{X}}_{0}}=\mathbf{A}$, which however has a great limitation that it requires $\left\| \mathbf{I}-\mathbf{A} \right\|<1$ to converge \supercite{17}. To overcome this limitation, we here recommend the following selection:
\begin{equation}\label{E12}
{{\mathbf{X}}_{0}}=\sqrt{\frac{\alpha }{\left\| \mathbf{A} \right\|}}\mathbf{A},
\end{equation}
where $\alpha >0$ is a given parameter. Then we have
\begin{equation}\label{E13}
{{\mathbf{Y}}_{0}}=\mathbf{I}-\frac{\alpha }{\left\| \mathbf{A} \right\|}\mathbf{A}
\end{equation}
According to the Theorem 4.15 in Ref.\supercite{26}, the convergence region of the iteration (\ref{E10}) depends on the following scalar iterative scheme
\begin{equation}\label{E14}
{{y}_{k+1}}=\frac{3}{4}y_{k}^{2}+\frac{1}{4}y_{k}^{3},\ {{y}_{0}}=1-\alpha a-\alpha b\text{i},
\end{equation}
where $a$ and $b$ are the real and imaginary parts of certain eigenvalue $z=a+b\text{i}$ of ${{\left\| \mathbf{A} \right\|}^{-1}}\mathbf{A}$, respectively. Due to $\rho \left( \mathbf{A} \right)\le \left\| \mathbf{A} \right\|$, we have $\left| a+b\text{i} \right|\le 1$. Obviously, when the eigenvalues of ${{\left\| \mathbf{A} \right\|}^{-1}}\mathbf{A}$ are all in
\begin{equation}\label{E15}
\left\{ z:\ z\in \mathbb{C}, \ \text{and} \ \left| z-\frac{1}{\alpha } \right|<\frac{1}{\alpha } \right\},
\end{equation}
$\left| {{y}_{0}} \right|<1$, and the scalar iteration (\ref{E14}) converges. Therefore, when the eigenvalues of ${{\left\| \mathbf{A} \right\|}^{-1}}\mathbf{A}$ are all real positive, taking $\alpha =2$ can insure the convergence of the iteration (\ref{E9}). In fact, the actual convergence region is larger than (\ref{E15}). Figure 1 displays the convergence regions when $\alpha $ takes different values. The plotting method is: select randomly ${{10}^{6}}$ initial points $y_{0}^{\left( i \right)}$ on the complex region $\left| z \right|<1$, and then substitute these initial points into (\ref{E14}) and iterate 200 steps. If for the $i$-th initial point $y_{0}^{\left( i \right)}$, $\left| y_{200}^{\left( i \right)} \right|\le {{10}^{-10}}$, $y_{0}^{\left( i \right)}$ is marked by red, or else $y_{0}^{\left( i \right)}$ is marked by blue. As can be seen in Figure 1, the value of $\alpha $ will significantly affect the size of convergence region. When $\alpha =0.5$, the right half region converges.
\begin{figure}\label{F1}
\centering
\includegraphics[width=145mm]{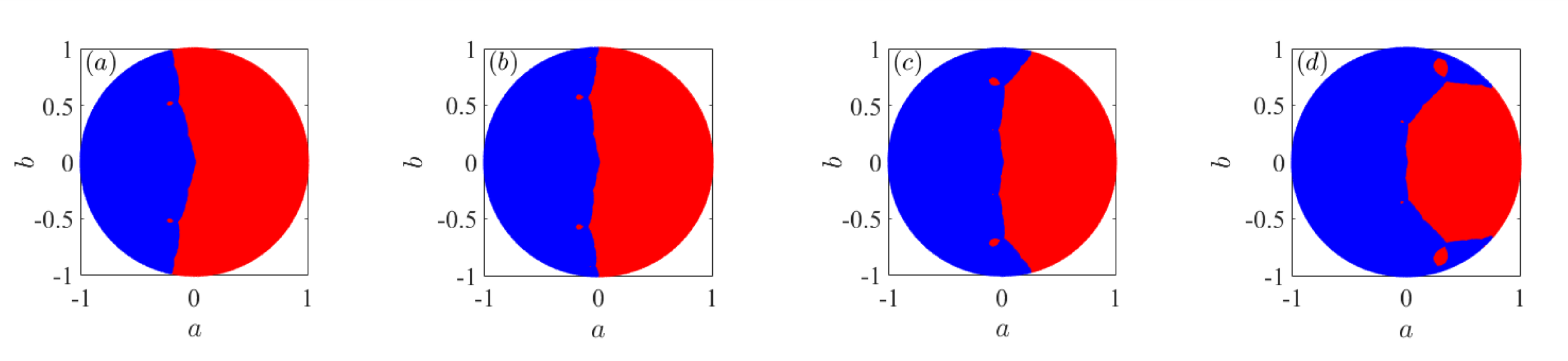}
\figcaption{Convergence region (red) for different $\alpha $: (a) $\alpha =0.1$; (b) $\alpha =0.5$; (c) $\alpha =1$; (d) $\alpha =2$.}
\end{figure}
\subsection{Stability analysis}\label{S2.2}
Although iteration (\ref{E9}) need to do only one matrix inversion, i.e., ${{\mathbf{A}}^{-1}}$, it is unstable as the iteration generated directly by the Newton method for (\ref{E7}), which will be proved below.

Adding a small error noted $\Delta {{\mathbf{X}}_{m}}$ at the $m$-th iteration, then we have ${{\mathbf{\hat{X}}}_{m}}={{\mathbf{X}}_{m}}+\Delta {{\mathbf{X}}_{m}}$, and 
\begin{equation}\label{E16}
{{\mathbf{\hat{X}}}_{m+1}}={{\mathbf{\hat{X}}}_{m}}\left( \mathbf{I}+\frac{1}{2}\left( \mathbf{I}-{{\mathbf{A}}^{-1}}\mathbf{\hat{X}}_{m}^{2} \right) \right)=\left( {{\mathbf{X}}_{m}}+\Delta {{\mathbf{X}}_{m}} \right)\left( \frac{3}{2}\mathbf{I}-\frac{1}{2}{{\mathbf{A}}^{-1}}{{\left( {{\mathbf{X}}_{m}}+\Delta {{\mathbf{X}}_{m}} \right)}^{2}} \right).
\end{equation}
Abandoning the high order terms of $\Delta {{\mathbf{X}}_{m}}$ gives
\begin{equation}\label{E17}
{{\mathbf{\hat{X}}}_{m+1}}={{\mathbf{X}}_{m+1}}-\frac{1}{2}{{\mathbf{A}}^{-1}}\left( \mathbf{X}_{m}^{2}\Delta {{\mathbf{X}}_{m}}+{{\mathbf{X}}_{m}}\Delta {{\mathbf{X}}_{m}}{{\mathbf{X}}_{m}} \right)+\Delta {{\mathbf{X}}_{m}}\left( \frac{3}{2}\mathbf{I}-\frac{1}{2}{{\mathbf{A}}^{-1}}\mathbf{X}_{m}^{2} \right).
\end{equation}
Noting that ${{\mathbf{\hat{X}}}_{m+1}}={{\mathbf{X}}_{m+1}}+\Delta {{\mathbf{X}}_{m+1}}$, we have
\begin{equation}\label{E18}
\Delta {{\mathbf{X}}_{m+1}}=-\frac{1}{2}{{\mathbf{A}}^{-1}}\left( \mathbf{X}_{m}^{2}\Delta {{\mathbf{X}}_{m}}+{{\mathbf{X}}_{m}}\Delta {{\mathbf{X}}_{m}}{{\mathbf{X}}_{m}} \right)+\Delta {{\mathbf{X}}_{m}}\left( \frac{3}{2}\mathbf{I}-\frac{1}{2}{{\mathbf{A}}^{-1}}\mathbf{X}_{m}^{2} \right).
\end{equation}
When $m\to \infty $, there will be $\left\| {{\mathbf{X}}_{m}} \right\|\to {{\mathbf{A}}^{0.5}}$, so we have ${{\mathbf{X}}_{m}}=\sqrt{\mathbf{A}}+{{\mathbf{\delta }}_{m}}$, where $\left\| {{\mathbf{\delta }}_{m}} \right\|\ll \left\| \sqrt{\mathbf{A}} \right\|$. Substituting ${{\mathbf{X}}_{m}}=\sqrt{\mathbf{A}}+{{\mathbf{\delta }}_{m}}$ into the above equation and ignoring the high-order small terms, there is
\begin{equation}\label{E19}
\Delta {{\mathbf{X}}_{m+1}}=\frac{1}{2}\Delta {{\mathbf{X}}_{m}}-\frac{1}{2}{{\mathbf{A}}^{-0.5}}\Delta {{\mathbf{X}}_{m}}{{\mathbf{A}}^{0.5}}.
\end{equation}

Applying the vec operator to (\ref{E19}), we have
\[\text{vec}\left( \Delta {{\mathbf{X}}_{m+1}} \right)=\mathbf{J}\text{vec}\left( \Delta {{\mathbf{X}}_{m}} \right),\ \ \mathbf{J}=\frac{1}{2}\left[ \mathbf{I}-{{\left( {{\mathbf{A}}^{\text{T}}} \right)}^{0.5}}\otimes {{\mathbf{A}}^{-0.5}} \right].\]
The stability condition of the above iterative scheme is $\rho \left( \mathbf{J} \right)\le 1$. However, according to the equation (6.24) in \supercite{26}, for Hermitian positive definite matrices, only if ${{\kappa }_{2}}\left( \mathbf{A} \right)<9$, $\rho \left( \mathbf{J} \right)\le 1$ will hold. So, the iteration (\ref{E9}) is conditional stable.

\subsection{The proposed iteration}\label{S2.3}

According to the above analysis, the iteration (\ref{E9}) cannot guarantee the stability. In this section we present a new stable iterative scheme. According to (\ref{E10}), ${{\mathbf{Y}}_{k+1}}$ can also be generated from ${{\mathbf{Y}}_{k}}$, then, the iterative scheme can be rewritten as
\begin{equation}\label{E20}
{{\mathbf{X}}_{k+1}}={{\mathbf{X}}_{k}}\left( \mathbf{I}+\frac{1}{2}{{\mathbf{Y}}_{k}} \right),\ \ \ {{\mathbf{Y}}_{k+1}}=\mathbf{Y}_{k}^{2}\left( \frac{3}{4}\mathbf{I}+\frac{1}{4}{{\mathbf{Y}}_{k}} \right).
\end{equation}
Obviously, if we do not take the round-off errors into consideration, the iteration schemes (\ref{E9}) and (\ref{E20}) are equivalent to each other, so the convergence analysis in subsection \ref{S2.2} is also applicable to (\ref{E20}). However, if the round off error is considered, we will find that the iteration (\ref{E20}) is numerically stable.

Due to the round-off error of computer, we gain the ${{\mathbf{\hat{X}}}_{m}}$ instead of ${{\mathbf{X}}_{m}}$ at the $m\text{-th}$ iteration, the relationship between the two is ${{\mathbf{\hat{X}}}_{m}}={{\mathbf{X}}_{m}}+\Delta {{\mathbf{X}}_{m}}$. Similarly, for ${{\mathbf{Y}}_{m}}$, we have ${{\mathbf{\hat{Y}}}_{m}}={{\mathbf{Y}}_{m}}+\Delta {{\mathbf{Y}}_{m}}$. So, there is:
\begin{equation}\label{E21}
\begin{aligned}
  {{{\mathbf{\hat{X}}}}_{m+1}}&= {{{\mathbf{\hat{X}}}}_{m}}\left( \mathbf{I}+0.5{{{\mathbf{\hat{Y}}}}_{m}} \right)=\left( {{\mathbf{X}}_{m}}+\Delta {{\mathbf{X}}_{m}} \right)\left( \mathbf{I}+0.5\left( {{\mathbf{Y}}_{m}}+\Delta {{\mathbf{Y}}_{m}} \right) \right) \\ 
 &=\text{ }{{\mathbf{X}}_{m+1}}+0.5{{\mathbf{X}}_{m}}\Delta {{\mathbf{Y}}_{m}}+\Delta {{\mathbf{X}}_{m}}\left( \mathbf{I}+0.5{{\mathbf{Y}}_{m}}+0.5\Delta {{\mathbf{Y}}_{m}} \right)  
\end{aligned}
\end{equation}
and
\begin{equation}\label{E22}
\begin{aligned}
 {{{\mathbf{\hat{Y}}}}_{m+1}} & =\mathbf{\hat{Y}}_{m}^{2}\left( \frac{3}{4}\mathbf{I}+\frac{1}{4}\mathbf{\hat{Y}}_{m}^{{}} \right)={{\left( {{\mathbf{Y}}_{m}}+\Delta {{\mathbf{Y}}_{m}} \right)}^{2}}\left( \frac{3}{4}\mathbf{I}+\frac{1}{4}\left( {{\mathbf{Y}}_{m}}+\Delta {{\mathbf{Y}}_{m}} \right) \right) \\ 
 & =\left( \mathbf{Y}_{m}^{2}+{{\mathbf{Y}}_{m}}\Delta {{\mathbf{Y}}_{m}}+\Delta {{\mathbf{Y}}_{m}}{{\mathbf{Y}}_{m}}+\Delta \mathbf{Y}_{m}^{2} \right)\left( \frac{3}{4}\mathbf{I}+\frac{1}{4}\left( {{\mathbf{Y}}_{m}}+\Delta {{\mathbf{Y}}_{m}} \right) \right). 
\end{aligned}
\end{equation}
Omit the second and higher order terms with respect to the round-off errors in (\ref{E21}) and (\ref{E22}), we have
\begin{align}
&\Delta {{\mathbf{X}}_{m+1}}=\Delta {{\mathbf{X}}_{m}}+0.5\Delta {{\mathbf{X}}_{m}}{{\mathbf{Y}}_{m}}+0.5{{\mathbf{X}}_{m}}\Delta {{\mathbf{Y}}_{m}},\label{E23}\\ 
&\Delta {{\mathbf{Y}}_{m+1}}=\frac{1}{4}\left( \mathbf{Y}_{m}^{2}\Delta {{\mathbf{Y}}_{m}}+{{\mathbf{Y}}_{m}}\Delta {{\mathbf{Y}}_{m}}{{\mathbf{Y}}_{m}}+\Delta {{\mathbf{Y}}_{m}}\mathbf{Y}_{m}^{2} \right)+\frac{3}{4}\left( {{\mathbf{Y}}_{m}}\Delta {{\mathbf{Y}}_{m}}+\Delta {{\mathbf{Y}}_{m}}{{\mathbf{Y}}_{m}} \right).\label{E24}
\end{align}
Take norms on both sides, we have
\begin{equation}\label{E25}
\left[ \begin{matrix}
   \left\| \Delta {{\mathbf{X}}_{m+1}} \right\|  \\
   \left\| \Delta {{\mathbf{Y}}_{m+1}} \right\|  \\
\end{matrix} \right]\le {{\mathbf{L}}_{m}}\left[ \begin{matrix}
   \left\| \Delta {{\mathbf{X}}_{m}} \right\|  \\
   \left\| \Delta {{\mathbf{Y}}_{m}} \right\|  \\
\end{matrix} \right],
\end{equation}
where
\begin{equation}\label{E26}
{{\mathbf{L}}_{m}}=\left[ \begin{matrix}
   1+0.5\left\| {{\mathbf{Y}}_{m}} \right\| & 0.5\left\| {{\mathbf{X}}_{m}} \right\|  \\
   0 & 1.5\left\| {{\mathbf{Y}}_{m}} \right\|+0.75{{\left\| {{\mathbf{Y}}_{m}} \right\|}^{2}}  \\
\end{matrix} \right].
\end{equation}
As $ m\to \infty $, $\left\| {{\mathbf{Y}}_{m}} \right\|\to 0$, and ${{\mathbf{L}}_{m}}$will become an idempotent matrix, which means the error will not grow infinitely. Hence the iteration (\ref{E20}) is stable.

It can be seen that the computation of ${{\mathbf{Y}}_{k+1}}$ is only depended on ${{\mathbf{Y}}_{k}}$, and no longer involve ${{\mathbf{A}}^{-1}}$. As ${{\mathbf{X}}_{0}}=\sqrt{0.5/\left\| \mathbf{A} \right\|}\mathbf{A}$, ${{\mathbf{Y}}_{0}}=\mathbf{I}-0.5\mathbf{A}/\left\| \mathbf{A} \right\|$ also have nothing to do with ${{\mathbf{A}}^{-1}}$. Hence the iterative scheme (\ref{E20}) combined with the initial matrix ${{\mathbf{X}}_{0}}=\sqrt{0.5/\left\| \mathbf{A} \right\|}\mathbf{A}$ is a stable iteration avoiding fully the matrix inversion. As matrix inversion costs much more time than matrix multiplication for large-scale matrices, avoiding matrix inversion can reduce computation time effectively.

\section{Sparsity of principal matrix square root}\label{S3}
Although the iterative algorithm (\ref{E20}) is stable and avoids computing the matrix inverse, it is still not applicable for the computation of principal square root of large-scale matrices due to the limitation of the computational cost and memory. However, it will be shown in this section that there are some sparse matrices whose principal square roots are still sparse, or nearly sparse, with the help of the bandwidth theory proposed in Ref.\supercite{22}. We first introduce the definitions of matrix bandwidth, real bandwidth and $\varepsilon $-bandwidth.
\begin{Dingli}[Definition 1]\label{D1}
For a sparse matrix $\mathbf{A}\in {{\mathbb{C}}^{n\times n}}$, the element of which is denoted ${{a}_{ij}}$, the bandwidth $l\left( \mathbf{A} \right)$ of the matrix $\mathbf{A}$ is defined by
\begin{equation}\label{E27}
l\left( \mathbf{A} \right)=\underset{\Omega \left( \mathbf{A} \right)}{\mathop{\max }}\,\left( j-i \right)+\underset{\Omega \left( \mathbf{A} \right)}{\mathop{\max }}\,\left( i-j \right),
\end{equation}
where $\Omega \left( \mathbf{A} \right):=\left\{ \left( i,j \right):{{a}_{ij}}\ne 0 \right\}$.
\end{Dingli}

\begin{Dingli}[Definition 2]\label{D2}
For a sparse matrix $\mathbf{A}\in {{\mathbb{C}}^{n\times n}}$, and any $\varepsilon >0$

(a)	the real bandwidth $\lambda \left( \mathbf{A} \right)$ is
\begin{equation}\label{E28}
\lambda \left( \mathbf{A} \right)=\underset{\mathbf{P}\in P}{\mathop{\min }}\,l\left( \mathbf{PA}{{\mathbf{P}}^{\text{T}}} \right).
\end{equation}
where $P$ is the set of all the permutation matrices;

(b)	the $\varepsilon $-bandwidth $\beta \left( \varepsilon ,\mathbf{A} \right)$ is
\begin{equation}\label{E29}
\beta \left( \varepsilon ,\mathbf{A} \right)=\underset{\mathbf{B}\in \Gamma \left( \varepsilon  \right)}{\mathop{\min }}\,\lambda \left( \mathbf{A}-\mathbf{B} \right),\text{  }\Gamma \left( \varepsilon  \right):=\left\{ \mathbf{B}:\mathbf{B}\in {{\mathbb{C}}^{n\times n}},\text{ and}\ \left\| \mathbf{B} \right\|\le \varepsilon \left\| \mathbf{A} \right\| \right\}.
\end{equation}
\end{Dingli}
For the real bandwidth and $\varepsilon $-bandwidth, we have the following lemma 3\supercite{22}. 
\begin{Dingli}[Lemma 3]\label{L3}
Let $\mathbf{A},\mathbf{B}\in {{\mathbb{C}}^{n\times n}}$ be sparse matrices, and ${{p}_{m}}\left( \mathbf{A} \right)=\sum\limits_{i=0}^{m}{{{c}_{i}}{{\mathbf{A}}^{i}}},{{c}_{i}}\in \mathbb{C}$,where $m$ is an arbitrary positive integer, then
 
(a) $\lambda \left( {{p}_{m}}\left( \mathbf{A} \right) \right)\le m\lambda \left( \mathbf{A} \right)$;

(b) $\beta \left( {{\varepsilon }_{p}},{{p}_{m}}\left( \mathbf{A} \right) \right)\le m\lambda \left( \mathbf{A}-\mathbf{B} \right)$, where ${{\varepsilon }_{p}}={{\left\| {{p}_{m}}\left( \mathbf{A} \right) \right\|}^{-1}}\left\| {{p}_{m}}\left( \mathbf{A} \right)-{{p}_{m}}\left( \mathbf{A}-\mathbf{B} \right) \right\|.$
\end{Dingli}

The sparsity of $\sqrt{\mathbf{A}}$ can be measured in terms of the real bandwidth or the $\varepsilon $-bandwidth. Actually, we have the following theorems:

\begin{Dingli}[Theorem 4]\label{T4}
For any sparse matrix $\mathbf{A}\in {{\mathbb{C}}^{n\times n}}$ without non-positive real eigenvalues and any $\varepsilon >0$, if ${{y}_{0}}=\left\| \mathbf{I}-0.5{{\left\| \mathbf{A} \right\|}^{-1}}\mathbf{A} \right\|<1,$ then
\begin{equation}\label{E30}
\beta \left( \varepsilon ,\sqrt{\mathbf{A}} \right)\le \nu \left( \varepsilon ,{{y}_{0}} \right)\lambda \left( \mathbf{A} \right),\ \nu \left( \varepsilon ,{{y}_{0}} \right)=\min \left\{ s,\ s\in \mathbb{N},\ \text{and},\ \sum\limits_{k=s}^{\infty }{\left| \frac{{{d}_{j}}}{j!} \right|y_{0}^{j}}\le \sqrt{0.5}\varepsilon  \right\}.
\end{equation}
\end{Dingli}

\begin{proof}
The Taylor series to $\sqrt{\mathbf{A}}$ can be expressed as
\begin{equation}\label{E31}
{{\mathbf{A}}^{1/2}}=\frac{\sqrt{\left\| \mathbf{A} \right\|}}{\sqrt{0.5}}\sqrt{\frac{0.5}{\left\| \mathbf{A} \right\|}\mathbf{A}}=\frac{\sqrt{\left\| \mathbf{A} \right\|}}{\sqrt{0.5}}\sqrt{\mathbf{I}-{{\mathbf{Y}}_{0}}}=\frac{\sqrt{\left\| \mathbf{A} \right\|}}{\sqrt{0.5}}\sum\limits_{j=0}^{\infty }{\frac{{{d}_{j}}}{j!}{{\left( -{{\mathbf{Y}}_{0}} \right)}^{j}}},\text{  }{{d}_{j}}=\prod\limits_{i=0}^{j-1}{\left( \frac{1}{2}-i \right)}.
\end{equation}
where ${{\mathbf{Y}}_{0}}=\mathbf{I}-0.5{{\left\| \mathbf{A} \right\|}^{-1}}\mathbf{A}$. So, there is a polynomial approximation to ${{\mathbf{A}}^{1/2}}$:
\begin{equation}\label{E32}
{{f}_{N}}\left( {{\mathbf{Y}}_{0}} \right)=\frac{\sqrt{\left\| \mathbf{A} \right\|}}{\sqrt{0.5}}\sum\limits_{j=0}^{N}{\frac{{{d}_{j}}}{j!}{{\left( -{{\mathbf{Y}}_{0}} \right)}^{j}}},
\end{equation}
whose real bandwidth satisfies $\lambda \left( {{f}_{N}}\left( {{Y}_{0}} \right) \right)\le N\lambda \left( \mathbf{A} \right)$, according to Lemma \ref{L3}(b). If $N$ satisfies
\begin{equation}\label{E33}
\frac{\left\| \frac{\sqrt{\left\| \mathbf{A} \right\|}}{\sqrt{0.5}}\sum\limits_{j=N+1}^{\infty }{\frac{{{d}_{j}}}{j!}{{\left( -{{\mathbf{Y}}_{0}} \right)}^{j}}} \right\|}{\left\| {{\mathbf{A}}^{1/2}} \right\|}\le \frac{1}{\sqrt{0.5}}\left\| \sum\limits_{j=N+1}^{\infty }{\frac{{{d}_{j}}}{j!}{{\left( -{{\mathbf{Y}}_{0}} \right)}^{j}}} \right\|\le \frac{1}{\sqrt{0.5}}\sum\limits_{j=N+1}^{\infty }{\left| \frac{{{d}_{j}}}{j!} \right|{{\left\| {{\mathbf{Y}}_{0}} \right\|}^{j}}}<\varepsilon, 
\end{equation}
there will be 
\begin{equation}\label{E34}
\frac{\left\| {{\mathbf{A}}^{1/2}}-{{f}_{N}}\left( {{\mathbf{Y}}_{0}} \right) \right\|}{\left\| {{\mathbf{A}}^{1/2}} \right\|}\le \varepsilon.
\end{equation}
Obviously, the minimum of $N$ satisfying (33) is $\nu \left( \varepsilon ,{{y}_{0}} \right)$. Then according to Lemma \ref{L3}(a), 
$\beta \left( \varepsilon ,\sqrt{\mathbf{A}} \right)\le \lambda \left( {{f}_{N}}\left( {{Y}_{0}} \right) \right)\le N\lambda \left( \mathbf{A} \right)\le \nu \left( \varepsilon ,{{y}_{0}} \right)\lambda \left( \mathbf{A} \right)$.
\end{proof}

\begin{Dingli}[Theorem 5]\label{T5}
For any sparse matrix $\mathbf{A}\in {{\mathbb{C}}^{n\times n}}$ with no non-positive real eigenvalues, the approximate ${{\mathbf{X}}_{j}}$ of the square root of $\mathbf{A}$ generated by the iteration (\ref{E20}) satisfies
\begin{equation}\label{E35}
\lambda \left( {{\mathbf{X}}_{j}} \right)\le \frac{{{3}^{j}}+1}{2}\lambda \left( \mathbf{A} \right),\ j=0,1,2,\cdots.
\end{equation}
\end{Dingli}

\begin{proof}
According to Lemma \ref{L3}(a) and considering that ${{\mathbf{X}}_{0}}=\sqrt{\frac{\alpha }{\left\| \mathbf{A} \right\|}}\mathbf{A}$, we have $\lambda \left( {{\mathbf{X}}_{0}} \right)\le \frac{{{3}^{0}}+1}{2}\lambda \left( \mathbf{A} \right)$. Supposing that $\lambda \left( {{\mathbf{X}}_{i}} \right)\le \frac{{{3}^{i}}+1}{2}\lambda \left( \mathbf{A} \right)$ holds for $i$, and considering the iteration (\ref{E20}), ${{\mathbf{X}}_{i+1}}={{\mathbf{X}}_{i}}+{{\mathbf{A}}^{-1}}\mathbf{X}_{i}^{3}$, ${{\mathbf{X}}_{i+1}}$ is an order-$\frac{{{3}^{i}}+1}{2}$ matrix polynomial of $\mathbf{A}$, then using Lemma 3.1(a), we have $\lambda \left( {{\mathbf{X}}_{i+1}} \right)\le \frac{{{3}^{i+1}}+1}{2}\lambda \left( \mathbf{A} \right)$. When $i=j$, $\lambda \left( {{\mathbf{X}}_{j}} \right)\le \frac{{{3}^{j}}+1}{2}\lambda \left( \mathbf{A} \right),\ j=0,1,2,\cdots.$
\end{proof}

According to Theorem 5, for any given relative error bound $\varepsilon >0$, the real bandwidth of the approximate matrix root calculated by (\ref{E20}) satisfies
\begin{equation}\label{E36}
\lambda \left( {{\mathbf{X}}_{\alpha }} \right)\le \frac{{{3}^{\alpha }}+1}{2}\lambda \left( \mathbf{A} \right),\alpha =\min \left\{ s:\frac{\left\| {{\mathbf{A}}^{1/2}}-{{\mathbf{X}}_{s}} \right\|}{\left\| {{\mathbf{A}}^{1/2}} \right\|}\le \varepsilon  \right\}.
\end{equation}
And according to Theorem 4, $\beta \left( \varepsilon ,\sqrt{\mathbf{A}} \right)\le \nu \left( \varepsilon ,{{y}_{0}} \right)\lambda \left( \mathbf{A} \right)$. Here we give a numerical example to compare the upper bound to the real bandwidth of the approximation matrix computed by the iteration (\ref{E20}) with that to the $\varepsilon $-bandwidth of ${{\mathbf{A}}^{1/2}}$ with $\varepsilon ={{10}^{-14}}$. The results are plotted in Figure 2. It can be seen that when $\varepsilon ={{10}^{-14}}$, we usually have
\begin{equation}\label{E37}
\nu \left( \varepsilon ,{{y}_{0}} \right)\ll\frac{{{3}^{\alpha }}+1}{2}.
\end{equation}
The numerical comparison shown in Figure 2 demonstrates that for the matrix computed by (\ref{E20}), there exists a sparser approximate matrix with the similar accuracy. This property is very helpful for the computation of large-scale sparse matrix square roots.
\begin{figure}\label{F2}
\centering
\includegraphics[width=70mm]{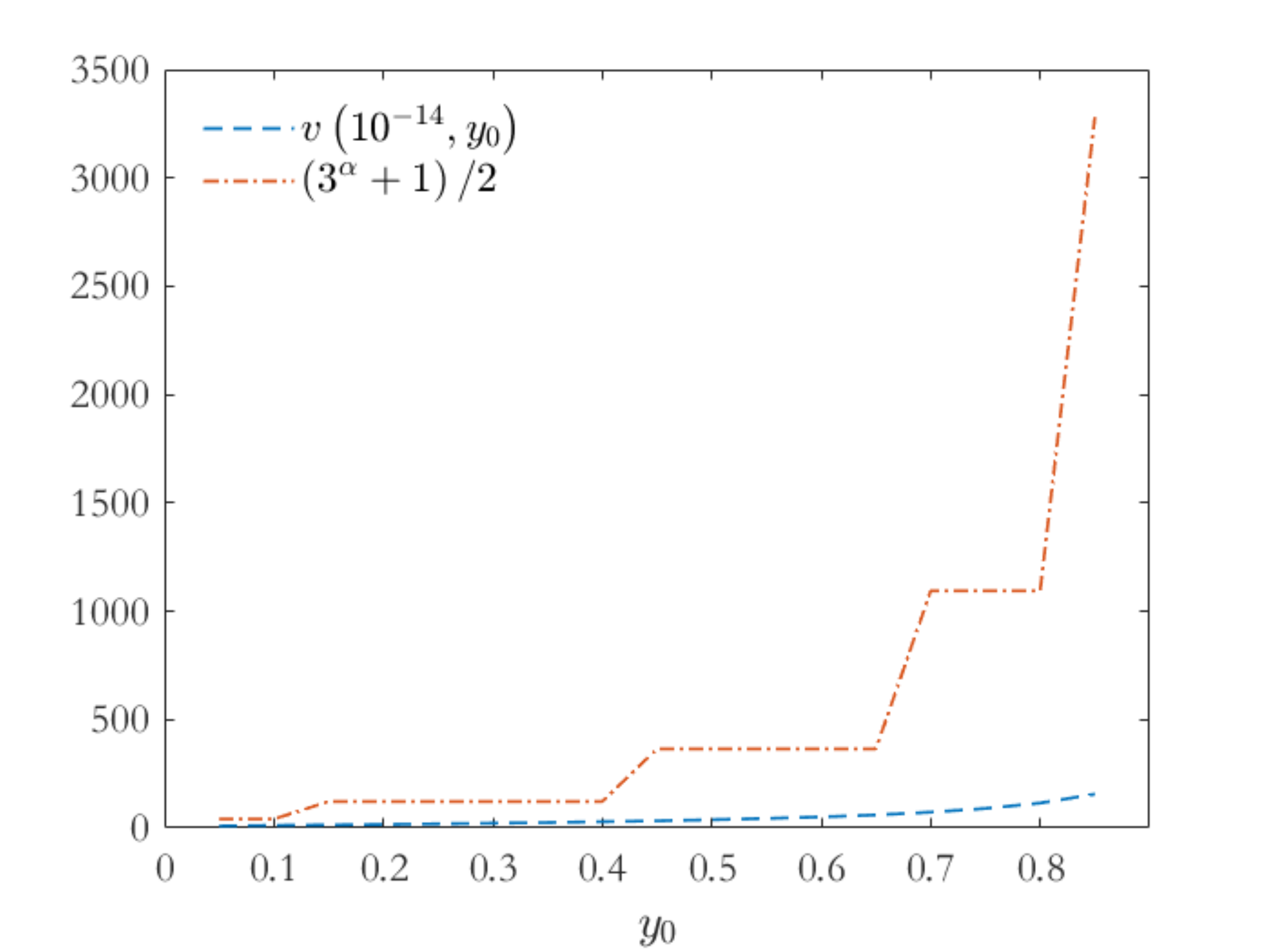}
\figcaption{Numerical comparison of $v\left( \varepsilon ,{{y}_{0}} \right)$ and $\left( {{3}^{\alpha }}+1 \right)/2$}
\end{figure}

\section{SIAI combined with filtering}\label{S4}
Obviously, if (\ref{E20}) is used directly, the results will be very dense due to multiple matrix multiplications, and the amount of computation and storage required are very large. Therefore, it is necessary to conduct in-depth research on how to obtain sparse approximate square roots. In this section, we will improve the iteration (\ref{E20}) by introducing an adaptive filtering computation to obtain a more efficient iterative method.
\subsection{Error analysis considering filtering}\label{S4.1}
Because the sparsity of  ${{\mathbf{X}}_{0}}=\sqrt{0.5/\left\| \mathbf{A} \right\|}\mathbf{A}$ and ${{\mathbf{Y}}_{0}}=\mathbf{I}-0.5\mathbf{A}/\left\| \mathbf{A} \right\|$ is the same as that of $\mathbf{A}$, so they are not needed to filter. Although both ${{\mathbf{X}}_{0}}$ and ${{\mathbf{Y}}_{0}}$ are sparse, the matrices ${{\mathbf{X}}_{1}}={{\mathbf{X}}_{0}}\left( \mathbf{I}+0.5{{\mathbf{Y}}_{0}} \right)$ and ${{\mathbf{Y}}_{1}}=\frac{3}{4}\mathbf{Y}_{0}^{2}+\frac{1}{4}\mathbf{Y}_{0}^{3}$ can become dense, because the real bandwidth is doubled due to the matrix multiplication. However, there exists a sparser matrix with the same precision as ${{\mathbf{X}}_{1}}$ according to the analysis in Section 3. Therefore, we apply the filtering technology to filter ${{\mathbf{X}}_{1}}$ and ignore the elements close to zeros, so as to obtain a sparser matrix ${{\mathbf{\hat{X}}}_{1}}={{\mathbf{X}}_{1}}-{{\mathbf{F}}_{1}}$, where ${{\mathbf{F}}_{1}}$ represents the matrix consisting of nearly zero elements dropped after ${{\mathbf{X}}_{1}}$ is filtered. Similarly, the filtering is still used to compute ${{\mathbf{Y}}_{1}}$, and its computation can be divided into
\begin{equation}\label{E38}
{{\mathbf{\hat{P}}}_{0}}=\mathbf{Y}_{0}^{2}-\Delta {{\mathbf{P}}_{0}}, \ {{\mathbf{\hat{Y}}}_{1}}={{\mathbf{\hat{P}}}_{0}}\left( \frac{3}{4}\mathbf{I}+\frac{1}{4}{{\mathbf{Y}}_{0}} \right)-{{\mathbf{E}}_{1}},
\end{equation}
where $\Delta {{\mathbf{P}}_{0}}$ and ${{\mathbf{E}}_{1}}$ represents the matrices consisting of nearly zero elements dropped after filtering $\mathbf{Y}_{0}^{2}$ and ${{\mathbf{\hat{P}}}_{0}}\left( \frac{3}{4}\mathbf{I}+\frac{1}{4}{{\mathbf{Y}}_{0}} \right)$, respectively. Subsequent calculations will be based on ${{\mathbf{\hat{X}}}_{1}}$ and ${{\mathbf{\hat{Y}}}_{1}}$. Generally, the new iterative format considering filtering (denoted by SIAI\_F) can be expressed as:
\begin{equation}\label{E39}
\left\{ \begin{aligned}
  & {{{\mathbf{\bar{X}}}}_{k+1}}={{{\mathbf{\hat{X}}}}_{k}}\left( \mathbf{I}+\frac{1}{2}{{{\mathbf{\hat{Y}}}}_{k}} \right),\ \ \ {{{\mathbf{\hat{X}}}}_{k+1}}={{{\mathbf{\bar{X}}}}_{k+1}}-{{\mathbf{F}}_{k+1}}, \\ 
 & {{\mathbf{P}}_{k}}=\mathbf{\hat{Y}}_{k}^{2},\ \ \ {{{\mathbf{\hat{P}}}}_{k}}={{\mathbf{P}}_{k}}-\Delta {{\mathbf{P}}_{k}}, \\ 
 & {{{\mathbf{\bar{Y}}}}_{k+1}}={{{\mathbf{\hat{P}}}}_{k}}\left( \frac{3}{4}\mathbf{I}+\frac{1}{4}{{{\mathbf{\hat{Y}}}}_{k}} \right),\ \ \ {{{\mathbf{\hat{Y}}}}_{k+1}}={{{\mathbf{\bar{Y}}}}_{k+1}}-{{\mathbf{E}}_{k+1}}. \\ 
\end{aligned} \right.
\end{equation}
Obviously, it is expected to apply appropriate constraints on ${{\mathbf{F}}_{k+1}}$, $\Delta {{\mathbf{P}}_{k}}$, and ${{\mathbf{E}}_{k+1}}$ to improve the computational efficiency as much as possible and ensure the accuracy at the same time. Let 
\begin{equation}\label{E40}
\Delta {{\mathbf{X}}_{k}}={{\mathbf{\hat{X}}}_{k}}-{{\mathbf{X}}_{k}}, \ \Delta {{\mathbf{Y}}_{k}}={{\mathbf{\hat{Y}}}_{k}}-{{\mathbf{Y}}_{k}}, \ k=0,1,2,\cdots.
\end{equation}
where ${{\mathbf{X}}_{k}}$ and ${{\mathbf{Y}}_{k}}$ are the matrices computed without filtering. They both meet the commutability with $\mathbf{A}$. Substituting (\ref{E40}) into the first equation of (\ref{E39}) and ignoring the higher-order small quantities, we have
\begin{equation}\label{E41}
\begin{aligned}
   {{{\mathbf{\bar{X}}}}_{k+1}}&={{{\mathbf{\hat{X}}}}_{k}}\left( \mathbf{I}+\frac{1}{2}{{{\mathbf{\hat{Y}}}}_{k}} \right)=\left( {{\mathbf{X}}_{k}}+\Delta {{\mathbf{X}}_{k}} \right)\left[ \mathbf{I}+\frac{1}{2}\left( {{\mathbf{Y}}_{k}}+\Delta {{\mathbf{Y}}_{k}} \right) \right] \\ 
 & ={{\mathbf{X}}_{k+1}}+\Delta {{\mathbf{X}}_{k}}\left( \mathbf{I}+\frac{1}{2}{{\mathbf{Y}}_{k}} \right)+\frac{1}{2}{{\mathbf{X}}_{k}}\Delta {{\mathbf{Y}}_{k}}.  
\end{aligned}
\end{equation}
Because ${{\mathbf{\hat{X}}}_{k}}={{\mathbf{\bar{X}}}_{k}}-{{\mathbf{F}}_{k}}$, (\ref{E41}) can be written as
\begin{equation}\label{E42}
{{\mathbf{\hat{X}}}_{k+1}}={{\mathbf{\bar{X}}}_{k+1}}-{{\mathbf{F}}_{k+1}}=\mathbf{X}_{k+1}^{{}}+\Delta {{\mathbf{X}}_{k}}\left( \mathbf{I}+\frac{1}{2}\mathbf{Y}_{k}^{{}} \right)+\frac{1}{2}{{\mathbf{X}}_{k}}\Delta {{\mathbf{Y}}_{k}}-{{\mathbf{F}}_{k+1}}.
\end{equation}
Combining the above equation with (\ref{E40}) yields
\begin{equation}\label{E43}
\Delta {{\mathbf{X}}_{k+1}}=\Delta {{\mathbf{X}}_{k}}\left( \mathbf{I}+\frac{1}{2}\mathbf{Y}_{k}^{{}} \right)+\frac{1}{2}{{\mathbf{X}}_{k}}\Delta {{\mathbf{Y}}_{k}}-{{\mathbf{F}}_{k+1}}.
\end{equation}

According to (\ref{E39}) and (\ref{E40}), ${{\mathbf{\hat{P}}}_{k}}$ can be expressed as
\begin{equation}\label{E44}
\begin{aligned}
  {{{\mathbf{\hat{P}}}}_{k}}& ={{\mathbf{P}}_{k}}-\Delta {{\mathbf{P}}_{k}}=\mathbf{\hat{Y}}_{k}^{2}-\Delta {{\mathbf{P}}_{k}}={{\left( {{\mathbf{Y}}_{k}}+\Delta {{\mathbf{Y}}_{k}} \right)}^{2}}-\Delta {{\mathbf{P}}_{k}} \\ 
 & =\mathbf{Y}_{k}^{2}+{{\mathbf{Y}}_{k}}\Delta {{\mathbf{Y}}_{k}}+\Delta {{\mathbf{Y}}_{k}}{{\mathbf{Y}}_{k}}-\Delta {{\mathbf{P}}_{k}}+\Delta \mathbf{Y}_{k}^{2}  
\end{aligned}
\end{equation}
Based on ${{\mathbf{\hat{P}}}_{k}}$, we can compute ${{\mathbf{\bar{Y}}}_{k+1}}$ by
\begin{equation}\label{E45}
{{\mathbf{\bar{Y}}}_{k+1}}=\mathbf{\hat{P}}_{k}^{{}}\left( \frac{3}{4}\mathbf{I}+\frac{1}{4}\mathbf{\hat{Y}}_{k}^{{}} \right)=\mathbf{\hat{P}}_{k}^{{}}\left[ \frac{3}{4}\mathbf{I}+\frac{1}{4}\left( \mathbf{Y}_{k}^{{}}+\Delta \mathbf{Y}_{k}^{{}} \right) \right].
\end{equation}
Substitute (\ref{E44}) into the above equation and ignore higher-order small quantities, we have
\begin{equation}\label{E46}
{{\mathbf{\bar{Y}}}_{k+1}}=\text{ }{{\mathbf{Y}}_{k+1}}+\frac{1}{4}\mathbf{Y}_{k}^{2}\Delta \mathbf{Y}_{k}^{{}}+\left( {{\mathbf{Y}}_{k}}\Delta {{\mathbf{Y}}_{k}}+\Delta {{\mathbf{Y}}_{k}}{{\mathbf{Y}}_{k}}-\Delta {{\mathbf{P}}_{k}} \right)\left( \frac{3}{4}\mathbf{I}+\frac{1}{4}\mathbf{Y}_{k}^{{}} \right).
\end{equation}
Because ${{\mathbf{\hat{Y}}}_{k}}={{\mathbf{\bar{Y}}}_{k}}-{{\mathbf{E}}_{k}}$, we have
\begin{equation}\label{E47}
{{\mathbf{\hat{Y}}}_{k+1}}={{\mathbf{\bar{Y}}}_{k+1}}-{{\mathbf{E}}_{k+1}}={{\mathbf{Y}}_{k+1}}+\frac{1}{4}\mathbf{Y}_{k}^{2}\Delta \mathbf{Y}_{k}^{{}}+\left( {{\mathbf{Y}}_{k}}\Delta {{\mathbf{Y}}_{k}}+\Delta {{\mathbf{Y}}_{k}}{{\mathbf{Y}}_{k}}-\Delta {{\mathbf{P}}_{k}} \right)\left( \frac{3}{4}\mathbf{I}+\frac{1}{4}\mathbf{Y}_{k}^{{}} \right)-{{\mathbf{E}}_{k+1}}.
\end{equation}
Combining the above equation with (\ref{E40}) gives
\begin{equation}\label{E48}
\Delta {{\mathbf{Y}}_{k+1}}=\frac{1}{4}\mathbf{Y}_{k}^{2}\Delta \mathbf{Y}_{k}^{{}}+\left( {{\mathbf{Y}}_{k}}\Delta {{\mathbf{Y}}_{k}}+\Delta {{\mathbf{Y}}_{k}}{{\mathbf{Y}}_{k}}-\Delta {{\mathbf{P}}_{k}} \right)\left( \frac{3}{4}\mathbf{I}+\frac{1}{4}\mathbf{Y}_{k}^{{}} \right)-{{\mathbf{E}}_{k+1}}.
\end{equation}
Due to ${{\mathbf{\hat{X}}}_{0}}={{\mathbf{X}}_{0}},{{\mathbf{\hat{Y}}}_{0}}={{\mathbf{Y}}_{0}}$, we have $\Delta {{\mathbf{X}}_{0}}=\mathbf{0}$, $\Delta {{\mathbf{Y}}_{0}}=\mathbf{0}$. Combining (\ref{E43}) and (\ref{E48}) yields
\begin{equation}\label{E49}
\left\{ \begin{aligned}
  & \Delta {{\mathbf{X}}_{k+1}}=\Delta {{\mathbf{X}}_{k}}\left( \mathbf{I}+0.5\mathbf{Y}_{k}^{{}} \right)+0.5{{\mathbf{X}}_{k}}\Delta {{\mathbf{Y}}_{k}}-{{\mathbf{F}}_{k+1}}, \\ 
 & \Delta {{\mathbf{Y}}_{k+1}}=0.25\mathbf{Y}_{k}^{2}\Delta \mathbf{Y}_{k}^{{}}+\left( {{\mathbf{Y}}_{k}}\Delta {{\mathbf{Y}}_{k}}+\Delta {{\mathbf{Y}}_{k}}{{\mathbf{Y}}_{k}}-\Delta {{\mathbf{P}}_{k}} \right)\left( 0.75\mathbf{I}+0.25\mathbf{Y}_{k}^{{}} \right)-{{\mathbf{E}}_{k+1}}, \\ 
& \Delta {{\mathbf{X}}_{0}}=\mathbf{0},\ \Delta {{\mathbf{Y}}_{0}}=\mathbf{0}.
\end{aligned} \right.
\end{equation}
From the above equation, we have
\begin{equation}\label{E50}
\left\{ \begin{aligned}
   \left\| \Delta {{\mathbf{X}}_{k}} \right\| \le &\left( \left\| \mathbf{I} \right\|+0.5\left\| \mathbf{Y}_{k-1}^{{}} \right\| \right)\left\| \Delta {{\mathbf{X}}_{k-1}} \right\|+0.5\left\| {{\mathbf{X}}_{k-1}} \right\|\left\| \Delta {{\mathbf{Y}}_{k-1}} \right\|+\left\| {{\mathbf{F}}_{k}} \right\|, \\ 
 \left\| \Delta {{\mathbf{Y}}_{k}} \right\|  \le &\left( 1.5\left\| \mathbf{I} \right\|+0.75\left\| \mathbf{Y}_{k-1}^{{}} \right\| \right)\left\| {{\mathbf{Y}}_{k-1}} \right\|\left\| \Delta {{\mathbf{Y}}_{k-1}} \right\| \\ & +\left( 0.75\left\| \mathbf{I} \right\|+0.25\left\| \mathbf{Y}_{k-1}^{{}} \right\| \right)\left\| \Delta {{\mathbf{P}}_{k-1}} \right\|+\left\| {{\mathbf{E}}_{k}} \right\|. \\ 
\end{aligned} \right.
\end{equation}
Rewrite (\ref{E50}) in the format of matrix,
\begin{equation}\label{E51}
{{\mathbf{U}}_{k}}\le \mathbf{T}_{k-1}^{{}}{{\mathbf{U}}_{k-1}}+{{\mathbf{B}}_{k-1}}\mathbf{V}_{k-1}^{{}},\  k=1,2,\cdots,
\end{equation}
where
\begin{equation}\label{E52}
\begin{aligned}
  & {{\mathbf{U}}_{k}}=\left[ \begin{matrix}
   \left\| \Delta {{\mathbf{X}}_{k}} \right\|  \\
   \left\| \Delta {{\mathbf{Y}}_{k}} \right\|  \\
\end{matrix} \right],\ \ \ {{\mathbf{T}}_{k-1}}=\left[ \begin{matrix}
   \left\| \mathbf{I} \right\|+0.5\left\| \mathbf{Y}_{k-1}^{{}} \right\| & 0.5\left\| {{\mathbf{X}}_{k-1}} \right\|  \\
   0 & \left( 1.5\left\| \mathbf{I} \right\|+0.75\left\| \mathbf{Y}_{k-1}^{{}} \right\| \right)\left\| \mathbf{Y}_{k-1}^{{}} \right\|  \\
\end{matrix} \right], \\ 
 & {{\mathbf{B}}_{k-1}}=\left[ \begin{matrix}
   0 & 0 & 1  \\
   0.75\left\| \mathbf{I} \right\|+0.25\left\| \mathbf{Y}_{k-1}^{{}} \right\| & 1 & 0  \\
\end{matrix} \right],\ \ \ \mathbf{V}_{k-1}^{{}}=\left[ \begin{matrix}
   \left\| \Delta {{\mathbf{P}}_{k-1}} \right\|  \\
   \left\| {{\mathbf{E}}_{k}} \right\|  \\
   \left\| {{\mathbf{F}}_{k}} \right\|  \\
\end{matrix} \right]. \\ 
\end{aligned}
\end{equation}
According to (\ref{E51}), if ${{\mathbf{U}}_{q}}$ is known, the error vector ${{\mathbf{U}}_{k}}$ will satisfy the following inequality: 
\begin{equation}\label{E53}
{{\mathbf{U}}_{k}}\le {{\mathbf{Q}}_{k,q}}{{\mathbf{U}}_{q}}+\sum\limits_{i=q+1}^{k}{{{\mathbf{Q}}_{k,i}}{{\mathbf{B}}_{i-1}}{{\mathbf{V}}_{i-1}}},\  {{\mathbf{Q}}_{k,i}}=\left\{ \begin{aligned}
   {{\mathbf{T}}_{k-1}}\cdots {{\mathbf{T}}_{i}}\text{ ,} \ &i=1,\ 2,\ \cdots ,\ k-1 \\ 
  \mathbf{I}\text{            ,} \ &i=k \\ 
\end{aligned} \right.
\end{equation}
which provides an evaluation for the accumulation of filtering errors.

\subsection{How to select adaptive filtering thresholds}\label{S4.2}
The upper bounds to $\left\| \Delta {{\mathbf{P}}_{k-1}} \right\|$, $\left\| {{\mathbf{E}}_{k}} \right\|$ and $\left\| {{\mathbf{F}}_{k}} \right\|$ are discussed here to ensure the accuracy of the computation results. After performing several iterations, the norm of the absolute error can be expressed as
\begin{equation}\label{E54}
\left\| {{\mathbf{A}}^{0.5}}-{{{\mathbf{\hat{X}}}}_{k}} \right\|=\left\| {{\mathbf{A}}^{0.5}}-{{\mathbf{X}}_{k}}+{{\mathbf{X}}_{k}}-{{{\mathbf{\hat{X}}}}_{k}} \right\|\le \left\| {{\mathbf{A}}^{0.5}}-{{\mathbf{X}}_{k}} \right\|+\left\| \Delta {{\mathbf{X}}_{k}} \right\|.
\end{equation}
Because ${{\mathbf{Y}}_{k}}=\mathbf{I}-{{\mathbf{A}}^{-1}}\mathbf{X}_{k}^{2}$, we have
\begin{equation}\label{E55}
{{\mathbf{A}}^{0.5}}={{\mathbf{X}}_{k}}{{\left( \mathbf{I}-{{\mathbf{Y}}_{k}} \right)}^{-0.5}}={{\mathbf{X}}_{k}}\left( \mathbf{I}+0.5{{\mathbf{Y}}_{k}}+O\left( \mathbf{Y}_{k}^{2} \right) \right).
\end{equation}
Substituting the above equation into (\ref{E54}) yields
\begin{equation}\label{E56}
\left\| {{\mathbf{A}}^{0.5}}-{{{\mathbf{\hat{X}}}}_{k}} \right\|\le \left\| 0.5{{\mathbf{X}}_{k}}{{\mathbf{Y}}_{k}} \right\|+\left\| \Delta {{\mathbf{X}}_{k}} \right\|+O\left( {{\left\| {{\mathbf{Y}}_{k}} \right\|}^{2}} \right).
\end{equation}
If the iteration converges after $M$ more iterations, we also have
\begin{equation}\label{E57}
\left\| {{\mathbf{A}}^{0.5}}-{{{\mathbf{\hat{X}}}}_{k+M}} \right\|\le \left\| 0.5{{\mathbf{X}}_{k+M}}{{\mathbf{Y}}_{k+M}} \right\|+\left\| \Delta {{\mathbf{X}}_{k+M}} \right\|+O\left( {{\left\| \mathbf{Y}_{k+M}^{{}} \right\|}^{2}} \right).
\end{equation}
Let $\mathbf{e}={{\left[ \begin{matrix}
   1 & 0  \\
\end{matrix} \right]}^{\text{T}}}$, according to (\ref{E53}), we have
\begin{equation}\label{E58}
\left\| \Delta {{\mathbf{X}}_{k+M}} \right\|\le {{\mathbf{e}}^{\text{T}}}{{\mathbf{U}}_{k+M}}={{\mathbf{e}}^{\text{T}}}\left( {{\mathbf{Q}}_{k+M,k}}{{\mathbf{U}}_{k}}+\sum\limits_{i=k+1}^{k+M}{{{\mathbf{Q}}_{k+M,i}}{{\mathbf{B}}_{i-1}}{{\mathbf{V}}_{i-1}}} \right).
\end{equation}
Substituting the above equation into (\ref{E57}) gives
\begin{equation}\label{E59}
\begin{aligned}
\left\| {{\mathbf{A}}^{0.5}}-{{{\mathbf{\hat{X}}}}_{k+M}} \right\|\le & \left\| 0.5{{\mathbf{X}}_{k+M}}{{\mathbf{Y}}_{k+M}} \right\| \\ 
&+{{\mathbf{e}}^{\text{T}}}\left( {{\mathbf{Q}}_{k+M,k}}{{\mathbf{U}}_{k}}+\sum\limits_{i=k+1}^{k+M}{{{\mathbf{Q}}_{k+M,i}}{{\mathbf{B}}_{i-1}}{{\mathbf{V}}_{i-1}}} \right)+O\left( {{\left\| \mathbf{Y}_{k+M}^{{}} \right\|}^{2}} \right).
\end{aligned}
\end{equation}
If the absolute error is required to be less than the given tolerance error ${{\varepsilon }_{\text{tol}}}$, then as long as
\begin{equation}\label{E60}
\left\| 0.5{{\mathbf{X}}_{k+M}}{{\mathbf{Y}}_{k+M}} \right\|+{{\mathbf{e}}^{\text{T}}}{{\mathbf{Q}}_{k+M,k}}{{\mathbf{U}}_{k}}+\sum\limits_{i=k+1}^{k+M}{{{\mathbf{e}}^{\text{T}}}{{\mathbf{Q}}_{k+M,i}}{{\mathbf{B}}_{i-1}}{{\mathbf{V}}_{i-1}}}+O\left( {{\left\| \mathbf{Y}_{k+M}^{{}} \right\|}^{2}} \right)\le {{\varepsilon }_{\text{tol}}},
\end{equation}
where ${{\mathbf{e}}^{\text{T}}}{{\mathbf{Q}}_{k+M,i}}{{\mathbf{B}}_{i-1}}{{\mathbf{V}}_{i-1}}$ denotes the error accumulation due to three filtrations in each iteration. If the error accumulation is required to be the same for each iteration, then we have
\begin{equation}\label{E61}
{{c}_{i,1}}\left\| \Delta {{\mathbf{P}}_{i-1}} \right\|+{{c}_{i,2}}\left\| {{\mathbf{E}}_{i}} \right\|+{{c}_{i,3}}\left\| {{\mathbf{F}}_{i}} \right\|\le \frac{{{\phi }_{k,M}}}{M}, \ {{\phi }_{k,M}}={{\varepsilon }_{\text{tol}}}-0.5\left\| {{\mathbf{X}}_{k+M}}{{\mathbf{Y}}_{k+M}} \right\|-{{\mathbf{e}}^{\text{T}}}{{\mathbf{Q}}_{k+M,k}}{{\mathbf{U}}_{k}},
\end{equation}
where \[{{c}_{i,1}}={{\mathbf{e}}^{\text{T}}}{{\mathbf{Q}}_{k+M,i}}{{\mathbf{B}}_{i-1}}{{\left[ \begin{matrix}
   1 & 0 & 0  \\
\end{matrix} \right]}^{\text{T}}},\]
\[{{c}_{i,2}}={{\mathbf{e}}^{\text{T}}}{{\mathbf{Q}}_{k+M,i}}{{\mathbf{B}}_{i-1}}{{\left[ \begin{matrix}
   0 & 1 & 0  \\
\end{matrix} \right]}^{\text{T}}},\]
\[{{c}_{i,3}}={{\mathbf{e}}^{\text{T}}}{{\mathbf{Q}}_{k+M,i}}{{\mathbf{B}}_{i-1}}{{\left[ \begin{matrix}
   0 & 0 & 1  \\
\end{matrix} \right]}^{\text{T}}}.\]
If we also use the idea that errors accumulate identically, one requires
\begin{equation}\label{E62}
\left\| \Delta {{\mathbf{P}}_{i-1}} \right\|\le \frac{{{\phi }_{k,M}}}{K{{c}_{i,1}}M},\ \ \ \left\| {{\mathbf{E}}_{i}} \right\|\le \frac{{{\phi }_{k,M}}}{K{{c}_{i,2}}M},\ \ \ \left\| {{\mathbf{F}}_{i}} \right\|\le \frac{{{\phi }_{k,M}}}{K{{c}_{i,3}}M},\ \ K=\sum\limits_{j=1}^{3}{\text{sign}\left( {{c}_{i,j}} \right)}.
\end{equation}
which is actually the adaptive filtering threshold and can be determined by simply estimating the number of residual iterations required from the current iteration to the final convergence. 

\subsection{How to estimate the number of residual iterations}\label{S4.3}
According to (\ref{E62}), $M$ must satisfy ${{\phi }_{k,M}}>0$, namely,
\begin{equation}\label{E63}
0.5\left\| {{\mathbf{X}}_{k+M}}{{\mathbf{Y}}_{k+M}} \right\|+{{\mathbf{e}}^{\text{T}}}{{\mathbf{Q}}_{k+M,k}}{{\mathbf{U}}_{k}}<{{\varepsilon }_{\text{tol}}}.
\end{equation}
For the convenience of further analysis, we denote ${{x}_{0}}=\left\| {{\mathbf{X}}_{k}} \right\|$ and ${{y}_{0}}=\left\| {{\mathbf{Y}}_{k}} \right\|$, and then do the following scalar iteration
\begin{equation}\label{E64}
\left\{ \begin{aligned}
  & {{x}_{m+1}}={{x}_{m}}+0.5{{x}_{m}}{{y}_{m}}, \\ 
 & {{y}_{m+1}}=y_{m}^{2}\left( 0.75+0.25{{y}_{m}} \right), \\ 
\end{aligned} \right.\ \ \ m=0,\ 1,\ \cdots ,\ M-1.
\end{equation}
In this study we use ${{x}_{M}}{{y}_{M}}$ to estimate $\left\| {{\mathbf{X}}_{k+M}}{{\mathbf{Y}}_{k+M}} \right\|$ and use ${{\mathbf{\tilde{Q}}}_{k+M,k}}$ to estimate ${{\mathbf{Q}}_{k+M,k}}$, where
\begin{equation}\label{E65}
{{\mathbf{\tilde{Q}}}_{k+M,k}}={{\mathbf{\tilde{T}}}_{k+M-1}}\cdots {{\mathbf{\tilde{T}}}_{k}},\ \ {{\mathbf{\tilde{T}}}_{k+m}}=\left[ \begin{matrix}
   1+0.5{{y}_{m}} & 0.5{{x}_{m}}  \\
   0 & \left( 1.5+0.75{{y}_{m}} \right){{y}_{m}}  \\
\end{matrix} \right].
\end{equation}
So, the number of residual iterations can be written as
\begin{equation}\label{E66}
M\left( k \right)=\min \left\{ m:\ \ m\in \mathbb{N}\ \ \text{and}\ \ 0.5{{x}_{m}}{{y}_{m}}+{{\mathbf{e}}^{\text{T}}}{{{\mathbf{\tilde{Q}}}}_{k+m,k}}{{\mathbf{U}}_{k}}\le {{\varepsilon }_{\text{tol}}} \right\}.
\end{equation}

Actually, since the accurate ${{\mathbf{X}}_{k}}$ and ${{\mathbf{Y}}_{k}}$ are unknown and the filtering threshold is usually small, we use ${{\mathbf{\hat{X}}}_{k}}$ and ${{\mathbf{\hat{Y}}}_{k}}$ instead, and then estimate ${{x}_{m}}$ and ${{y}_{m}}$ based on (\ref{E64}). In theory, the total number of iterations required for the convergence is a constant, and hence $k+M\left( k \right)$ must be a constant. With the increase of $k$, $M\left( k \right)$ will decrease. However, it is very conservative to use ${{x}_{m}}$ and ${{y}_{m}}$ to estimate $\left\| {{\mathbf{X}}_{k+M}}{{\mathbf{Y}}_{k+M}} \right\|$ and $M\left( k \right)$, it cannot guarantee that $k+M\left( k \right)$ is a constant. In the first few steps of the iteration, we use ${{x}_{m}}$ and ${{y}_{m}}$ to estimate $\left\| {{\mathbf{X}}_{k+M}}{{\mathbf{Y}}_{k+M}} \right\|$, which necessary leads to a large difference between $\left\| {{\mathbf{X}}_{k+M}}{{\mathbf{Y}}_{k+M}} \right\|$ and ${{x}_{m}}{{y}_{m}}$. But with the iteration going on, this difference will change to be smaller and smaller, and $k+M\left( k \right)$ will be close to a constant. When $k+M\left( k \right)$ becomes a constant, we believe that the current estimation of $M$ is accurate, and the iteration will converge after $M$ iterations.

According to Section \ref{S2.1}, as long as $\rho \left( {{\mathbf{Y}}_{0}} \right)<1$, the iteration will converge, but sometimes $\left\| {{\mathbf{Y}}_{0}} \right\|$ may be greater than 1, which makes $M$ have no solution when using the estimation (\ref{E66}), and the adaptive filtering thresholds (\ref{E62}) cannot be used. In this case, a fixed threshold is recommended as a supplement, namely, when $\left\| {{\mathbf{Y}}_{k}} \right\|>0.96$, it is required that
\begin{equation}\label{E67}
\left\| \Delta {{\mathbf{P}}_{k-1}} \right\|\le 0.01{{\varepsilon }_{\text{tol}}},\ \ \ \left\| {{\mathbf{E}}_{k}} \right\|\le 0.01{{\varepsilon }_{\text{tol}}},\ \ \ \left\| {{\mathbf{F}}_{k}} \right\|\le 0.01{{\varepsilon }_{\text{tol}}}.
\end{equation}

\subsection{Adaptive algorithm using relative error bound}\label{S4.4}
In the actual computation, the relative error bound rather than the absolute error bound  is often used as the tolerance error. In the previous section, we use the absolute error bound as the tolerance error, and gives the corresponding adaptive filtering threshold. If the relative error bound is selected as the tolerance error ${{\varepsilon }_{\text{tol}}}$, then according to (\ref{E60}) it is required that
\begin{equation}\label{E68}
\left\| 0.5{{\mathbf{X}}_{k+M}}{{\mathbf{Y}}_{k+M}} \right\|+{{\mathbf{e}}^{\text{T}}}{{\mathbf{Q}}_{k+M,k}}{{\mathbf{U}}_{k}}+\sum\limits_{i=k+1}^{k+M}{{{\mathbf{e}}^{\text{T}}}{{\mathbf{Q}}_{k+M,i}}{{\mathbf{B}}_{i-1}}{{\mathbf{V}}_{i-1}}}+O\left( {{\left\| \mathbf{Y}_{k+M}^{{}} \right\|}^{2}} \right)\le \left\| \sqrt{\mathbf{A}} \right\|{{\varepsilon }_{\text{tol}}}.
\end{equation}
where $\left\| \sqrt{\mathbf{A}} \right\|$ can be evaluated in terms of the following lemma.
\begin{Dingli}[Lemma 7]\label{L7}
For the matrices ${{\mathbf{X}}_{k}}$ and ${{\mathbf{Y}}_{k}}$ in the iterative scheme (\ref{E20}), if $\left\| {{\mathbf{Y}}_{k}} \right\|\le 1$, then
\begin{equation}\label{E69}
\left\| \sqrt{\mathbf{A}} \right\|\ge {{a}_{k}}:=\max \left( \sqrt{\left\| \mathbf{A} \right\|},\ \ \frac{\left\| {{\mathbf{X}}_{k}} \right\|}{2-\sqrt{1-\left\| {{\mathbf{Y}}_{k}} \right\|}} \right).
\end{equation}
\end{Dingli}
\begin{proof}
According to ${{\mathbf{Y}}_{k}}=\mathbf{I}-{{\mathbf{A}}^{-1}}\mathbf{X}_{k}^{2}$, we have
\[{{\mathbf{X}}_{k}}=\sqrt{\mathbf{A}}\sqrt{\mathbf{I}-{{\mathbf{Y}}_{k}}}=\sqrt{\mathbf{A}}\left( \mathbf{I}-\frac{1}{2}{{\mathbf{Y}}_{k}}+\frac{1}{8}\mathbf{Y}_{k}^{2}-\cdots  \right),\]
which means
\begin{equation}\label{E70}
\left\| {{\mathbf{X}}_{k}} \right\|\le \left\| \sqrt{\mathbf{A}} \right\|\left( 1+\frac{1}{2}\left\| {{\mathbf{Y}}_{k}} \right\|+\frac{1}{8}{{\left\| {{\mathbf{Y}}_{k}} \right\|}^{2}}+... \right)=\left\| \sqrt{\mathbf{A}} \right\|\left( 2-\sqrt{1-\left\| {{\mathbf{Y}}_{k}} \right\|} \right).
\end{equation}
As $\left\| \mathbf{A} \right\|=\left\| {{\left( \sqrt{\mathbf{A}} \right)}^{2}} \right\|\le {{\left\| \left( \sqrt{\mathbf{A}} \right) \right\|}^{2}}$, there will be $\left\| \sqrt{\mathbf{A}} \right\|\ge \sqrt{\left\| \mathbf{A} \right\|}$ and the combining of which with (\ref{E70}) completes the proof.
\end{proof}

In the actual computation, ${{a}_{k}}$ is evaluated by using
\begin{equation}\label{E71}
{{a}_{k}}=\max \left( \sqrt{\left\| \mathbf{A} \right\|},\ \ \frac{\left\| {{{\mathbf{\hat{X}}}}_{k}} \right\|}{2-\sqrt{1-\left\| {{{\mathbf{\hat{Y}}}}_{k}} \right\|}} \right).
\end{equation}
Using Lemma 7, (\ref{E68}) can be rewritten as
\begin{equation}\label{E72}
\left\| 0.5{{\mathbf{X}}_{k+M}}{{\mathbf{Y}}_{k+M}} \right\|+{{\mathbf{e}}^{\text{T}}}{{\mathbf{Q}}_{k+M,k}}{{\mathbf{U}}_{k}}+\sum\limits_{i=k+1}^{k+M}{{{\mathbf{e}}^{\text{T}}}{{\mathbf{Q}}_{k+M,i}}{{\mathbf{B}}_{i-1}}{{\mathbf{V}}_{i-1}}}+O\left( {{\left\| \mathbf{Y}_{k+M}^{{}} \right\|}^{2}} \right)\le {{a}_{k}}{{\varepsilon }_{\text{tol}}}.
\end{equation}
Through the same derivation processes proposed in subsection 4.2, the upper bounds used for the filtering are 
\begin{equation}\label{E73}
\left\| \Delta {{\mathbf{P}}_{i-1}} \right\|\le \frac{{{{\tilde{\phi }}}_{k,M}}}{K{{c}_{i,1}}M},\ \ \ \left\| {{\mathbf{E}}_{i}} \right\|\le \frac{{{{\tilde{\phi }}}_{k,M}}}{K{{c}_{i,2}}M},\ \ \ \left\| {{\mathbf{F}}_{i}} \right\|\le \frac{{{{\tilde{\phi }}}_{k,M}}}{K{{c}_{i,3}}M},\ \ K=\sum\limits_{j=1}^{3}{\text{sign}\left( {{c}_{i,j}} \right)}.
\end{equation}
where
\[{{\tilde{\phi }}_{k,M}}={{a}_{k}}{{\varepsilon }_{\text{tol}}}-0.5\left\| {{\mathbf{X}}_{k+M}}{{\mathbf{Y}}_{k+M}} \right\|-{{\mathbf{e}}^{\text{T}}}{{\mathbf{Q}}_{k+M,k}}{{\mathbf{U}}_{k}}.\]

\subsection{Specific algorithm of the proposed method}\label{S4.5}
\begin{algorithm}
\caption{Determination of the adaptive thresholds $\\$ [${{p}_{\text{k,ft}}}$, ${{e}_{\text{k,ft}}}$, ${{f}_{\text{k,ft}}}$] = adaptive\_thre[$\left\| {{{\mathbf{\hat{X}}}}_{k-1}} \right\|$, $\left\| {{{\mathbf{\hat{Y}}}}_{k-1}} \right\|$, ${{\mathbf{U}}_{k-1}}$, ${{\mathbf{B}}_{k-1}}$, $C$, ${{\varepsilon }_{\text{tol}}}$, ${{a}_{k-1}}$]}
\begin{algorithmic}[1]\label{A1}
\IF {$\left\| {{{\mathbf{\hat{Y}}}}_{k-1}} \right\|>0.96$}
\RETURN ${{p}_{\text{k,ft}}}=C{{\varepsilon }_{\text{tol}}}$; ${{e}_{\text{k,ft}}}=C{{\varepsilon }_{\text{tol}}}$; ${{f}_{\text{k,ft}}}=C{{\varepsilon }_{\text{tol}}}$;
\ELSE
\STATE $M=0$; ${{\mathbf{\tilde{Q}}}_{k+M-1,k-1}}={{\mathbf{I}}_{2\times 2}}$; ${{x}_{0}}=\left\| {{{\mathbf{\hat{X}}}}_{k-1}}\right\|$; ${{y}_{0}}=\left\| {{{\mathbf{\hat{Y}}}}_{k-1}} \right\|$;
\STATE ${{\phi }_{k-1,0}}={{a}_{k-1}}{{\varepsilon }_{\text{tol}}}-0.5{{x}_{0}}{{y}_{0}}-\left[ \begin{matrix}
   1 & 0  \\
\end{matrix} \right]{{\mathbf{U}}_{k-1}}$;
\WHILE {${{\phi }_{k-1,M}}\le 0$}
\STATE $M=M+1$;
\STATE Compute ${{\mathbf{\tilde{T}}}_{k+M-1}}$ according to (\ref{E65});
\STATE ${{\mathbf{\tilde{Q}}}_{k+M-1,k-1}}={{\mathbf{\tilde{T}}}_{k+M-1}}{{\mathbf{\tilde{Q}}}_{k+M-2,k-1}}$;
\IF {$M=1$}
\STATE ${{\mathbf{\tilde{Q}}}_{k+M-1,k}}={{\mathbf{I}}_{2\times 2}}$;
\ELSE
\STATE ${{\mathbf{\tilde{Q}}}_{k+M-1,k}}={{\mathbf{\tilde{T}}}_{k+M-1,k}}{{\mathbf{\tilde{Q}}}_{k+M-1,k}}$;
\ENDIF
\STATE ${{x}_{M}}={{x}_{M-1}}+0.5{{x}_{M-1}}{{y}_{M-1}}$; ${{y}_{M}}=y_{M-1}^{2}\left( 0.75+0.25{{y}_{M-1}} \right)$;
\STATE ${{\phi }_{k-1,M}}={{a}_{k-1}}{{\varepsilon }_{\text{tol}}}-0.5{{x}_{M}}{{y}_{M}}-\left[ \begin{matrix}
   1 & 0  \\
\end{matrix} \right]{{\mathbf{\tilde{Q}}}_{k+M-1,k-1}}{{\mathbf{U}}_{k-1}}$;
\ENDWHILE
\STATE Compute ${{c}_{k,1}}$; ${{c}_{k,2}}$; ${{c}_{k,3}}$ according to (\ref{E61});
\RETURN ${{p}_{\text{k,ft}}}=\frac{{{\phi }_{k-1,M}}}{K{{c}_{k,1}}M}$; ${{e}_{\text{k,ft}}}=\frac{{{\phi }_{k-1,M}}}{K{{c}_{k,2}}M}$; ${{f}_{\text{k,ft}}}=\ \frac{{{\phi }_{k-1,M}}}{K{{c}_{k,3}}M}$; $K=\sum\limits_{j=1}^{3}{\text{sign}\left( {{c}_{k,j}} \right)}$;
\ENDIF
\end{algorithmic}
\end{algorithm}

\begin{algorithm}[H]
\caption{This algorithm computes the $\sqrt{\mathbf{A}}$ for the given matrix $\mathbf{A}$. The (absolute/relative) error of the result will be smaller than the given tolerance error ${{\varepsilon }_{\text{tol}}}$.}
\begin{algorithmic}[1]\label{A2}
\STATE Let ${{\mathbf{\hat{X}}}_{0}}=\sqrt{0.5/\left\| \mathbf{A} \right\|}\mathbf{A}$; ${{\mathbf{\hat{Y}}}_{0}}=\mathbf{I}-0.5\mathbf{A}/\left\| \mathbf{A} \right\|$; ${{\mathbf{\hat{S}}}_{0}}=0.5{{\mathbf{\hat{X}}}_{0}}{{\mathbf{\hat{Y}}}_{0}}$;
\STATE ${{y}_{0}}=0$; ${{x}_{0}}=0$; ${{\mathbf{U}}_{0}}={{\left[ \begin{matrix}
   {{x}_{0}} & {{y}_{0}}  \\
\end{matrix} \right]}^{\text{T}}}$; ${{\mathbf{B}}_{0}}=\left[ \begin{matrix}
   0 & 0 & 1  \\
   0.75+\left\| {{{\mathbf{\hat{Y}}}}_{0}} \right\| & 1 & 0  \\
\end{matrix} \right]$; $C=0.01$;
\STATE $k=0$; ${{a}_{k}}=1$ (absolute error bound) or ${{a}_{k}}=\sqrt{\left\| \mathbf{A} \right\|}$ (relative error bound);
\WHILE {$\left\| {{{\mathbf{\hat{S}}}}_{k}} \right\|+{{x}_{k}}>{{a}_{k}}{{\varepsilon }_{\text{tol}}}$}
\STATE $k=k+1;$ 
\STATE  [${{p}_{\text{k,ft}}}$,${{e}_{\text{k,ft}}}$,${{f}_{\text{k,ft}}}$] = adaptive\_thre [$\left\| {{{\mathbf{\hat{X}}}}_{k-1}} \right\|$,$\left\| {{{\mathbf{\hat{Y}}}}_{k-1}} \right\|$,${{\mathbf{U}}_{k-1}}$,${{\mathbf{B}}_{k-1}}$,$C$,${{\varepsilon }_{\text{tol}}}$,${{a}_{k-1}}$];
\STATE $\text{Filt out }{{\mathbf{P}}_{k-1}}\text{ to get }{{\mathbf{\hat{P}}}_{k-1}} \ \text{which}\ \text{satisfies}\left\| \Delta {{\mathbf{P}}_{k-1}} \right\|=\left\| {{{\mathbf{\hat{P}}}}_{k-1}}-{{\mathbf{P}}_{k-1}} \right\|\le {{p}_{\text{k,ft}}}$;
\STATE ${{\mathbf{\bar{Y}}}_{k}}={{\mathbf{\hat{P}}}_{k-1}}\left( 0.75\mathbf{I}+0.25\mathbf{\hat{Y}}_{k-1}^{{}} \right)$;
\STATE $\text{Filt out }{{\mathbf{\bar{Y}}}_{k}}\text{ to get }{{\mathbf{\hat{Y}}}_{k}}\ \text{which}\ \text{satisfies}\ \left\| {{\mathbf{E}}_{k}} \right\|=\left\| {{{\mathbf{\hat{Y}}}}_{k}}-{{{\mathbf{\bar{Y}}}}_{k}} \right\|\le {{e}_{\text{k,ft}}}$;
\STATE ${{\mathbf{\bar{X}}}_{k}}={{\mathbf{\hat{X}}}_{k-1}}+{{\mathbf{\hat{S}}}_{k-1}}$;
\STATE $\text{Filt out }{{\mathbf{\bar{X}}}_{k}}\text{ to get }{{\mathbf{\hat{X}}}_{k}}\text{which}\ \text{satisfies}\left\| {{\mathbf{F}}_{k}} \right\|=\left\| {{{\mathbf{\hat{X}}}}_{k}}-{{{\mathbf{\bar{X}}}}_{k}} \right\|\le {{f}_{\text{k,ft}}}$;
\STATE ${{\mathbf{\hat{S}}}_{k}}=0.5{{\mathbf{\hat{X}}}_{k}}{{\mathbf{\hat{Y}}}_{k}}$;
\STATE Compute ${{\mathbf{T}}_{k-1}}$, ${{\mathbf{B}}_{k-1}}$ and ${{\mathbf{V}}_{k-1}}$ according to (\ref{E52});
\STATE ${{\mathbf{U}}_{k}}={{\mathbf{T}}_{k-1}}{{\mathbf{U}}_{k-1}}+{{\mathbf{B}}_{k-1}}{{\mathbf{V}}_{k-1}}$;
\STATE ${{x}_{k}}=\left[ \begin{matrix}
   1 & 0  \\
\end{matrix} \right]{{\mathbf{U}}_{k}}$;
\IF {$\left\| {{{\mathbf{\hat{Y}}}}_{k}} \right\|<1$}
\STATE ${{a}_{k}}=\max \left( \sqrt{\left\| \mathbf{A} \right\|},\ \ {\left\| {{{\mathbf{\hat{X}}}}_{k}} \right\|}/{\left( 2-\sqrt{1-\left\| {{{\mathbf{\hat{Y}}}}_{k}} \right\|} \right)}\; \right)$;
\ENDIF
\ENDWHILE
\RETURN ${{\mathbf{A}}^{1/2}}\approx {{\mathbf{\hat{X}}}_{k}}$;
\end{algorithmic} 
\end{algorithm}
The SIAI\_F can be summarized as the Algorithm 1 and Algorithm 2. In lines 7, 9 and 11 of Algorithm \ref{A2}, we use the filtering algorithm proposed in Ref.\supercite{22} to filter the matrix to obtain the sparser matrix. If it is required to use the absolute error bound, the lines 16-18 in Algorithm \ref{A2} should be deleted. 

\section{Numerical experiments} \label{S5}
To test the accuracy of the proposed algorithm, we define the following relative error:
\begin{equation}\label{E74}
er=\frac{\left\| {{\mathbf{X}}^{2}}-\mathbf{A} \right\|}{\left\| \mathbf{A} \right\|}.
\end{equation}

Due to the difficulty of the evaluation of the 2-norm of a large-scale matrix, 1-norm is applied here. We performed the experiments by using the computer with Microsoft Windows 10 21H1, Intel(R) Core (TM) i7-10750H CPU @ 2.60GHz, and 31.6GB of RAM. And the MATLAB version is MATLAB R2021b. In this section, we note the iteration (20) as SIAI and note the Algorithm 4.2 as SIAI\_F for convenience. The code of SIAI\_F, as well as the matrices used in the experiment 3, are uploaded to website https://www.rocewea.com/8.html for interested readers to download.

\subsection{Performance of the SIAI}\label{S5.1}
As the first experiment, we compare the proposed SIAI with the original algorithm without filtering, i.e., the iterative scheme (\ref{E9}), to test the effect of the matrix inversion on the computational efficiency. In addition, the DB method, the IN method and the sqrtm function$^{[2,27\textrm{--}28]}$ in MATLAB are also performed. Both the DB and the IN need matrix inversions. All the four methods are performed with the format of the full matrix. The considered matrix $\mathbf{A}$, which is related to the discretization of partial differential equations, is
\begin{equation}\label{E75}
\mathbf{A}=\left[ \begin{matrix}
   1-2\lambda  & \lambda  & {} & {} & {}  \\
   \lambda  & 1-2\lambda  & \lambda  & {} & {}  \\
   {} & \ddots  & \ddots  & \ddots  & {}  \\
   {} & {} & \lambda  & 1-2\lambda  & \lambda   \\
   {} & {} & {} & \lambda  & 1-2\lambda   \\
\end{matrix} \right].
\end{equation}
where $\lambda =-1$, and $n=500,1000,1500,2000$. The computational time and errors of different methods are listed in Table \ref{Tab1}.

\begin{center}
\abovecaptionskip 0pt \belowcaptionskip 1pt
\renewcommand{\arraystretch}{1.2}
{\scriptsize \tabcaption{Comparisons of different algorithms}\label{Tab1}
\noindent\begin{tabular*}{\textwidth}{@{\extracolsep{\fill}}@{~~}ccccccccc}
\toprule
\raisebox{-2.00ex}[0cm][0cm]{Method} &
\multicolumn{2}{c}{$n=500$} & \multicolumn{2}{c}{$n=1000$}& \multicolumn{2}{c}{$n=1500$} & \multicolumn{2}{c}{$n=2000$}\\
\cline{2-9}
& $er$ & Time(s) & $er$ & Time(s) & $er$ & Time(s) & $er$ & Time(s)\\
\midrule
SIAI & 1.42e-15 & 0.149 & 1.42e-15 & 0.848 & 1.42e-15 & 2.172 & 1.42e-15 & 4.331\\
Iteration (\ref{E9})& 9.84e-16 & 0.264 & 9.84e-16 & 2.918 & 1.29e-15 & 7.593 & 9.96e-16 & 13.697\\
DB & 1.39e-15 & 0.136 & 1.39e-15 & 1.731 & 1.39e-15 & 3.209 & 1.39e-15 & 5.855\\
IN & 5.52e-16 & 0.159 & 5.52e-16 & 1.983 & 5.52e-16 & 4.974 & 5.52e-16 & 9.225\\
sqrtm & 4.41e-14 & 0.008 & 8.31e-14 & 0.037 & 1.13e-13 & 0.100 & 1.45e-13 & 0.207\\
\bottomrule
\end{tabular*}\small }\\[4mm]
\end{center}

As is shown in Table \ref{Tab1}, the SIAI is more efficient than the original algorithm (Iteration (\ref{E9})), the DB, and the IN, but the accuracy of these four methods are the same. The improvement of the computational efficiency is due to avoiding the matrix inversion. Among the five methods, the sqrtm has the best calculation efficiency, but its accuracy is also the worst. This is because sqrtm is based on Schur decomposition, which is essentially a direct solution. Obviously, SIAI is a better choice when high-precision computations are needed. It must be pointed out that in this experiment the maximum dimension is only 2000, hence the full matrix instead of sparse matrix is used in the computation. In the following experiments, we will further study the computation of sparse matrices.

\subsection{Changes of the sparse matrix bandwidth in iteration}\label{S5.2}
The second experiment focus on the change of the bandwidths of the matrices involved in the proposed method by using the same matrix defined by (\ref{E75}) with $\lambda =-1$ and $n=10000$. The SIAI, and the SIAI\_F are applied to the computation of this matrix. The SIAI\_F used here is the relative error version of algorithm \ref{A2} with the tolerance relative error ${{\varepsilon }_{\text{tol}}}={{10}^{-13}}$. 

Both the SIAI and the SIAI\_F require 7 iterations under the given tolerance error. In Table \ref{Tab2}, the real bandwidths $\lambda \left( {{\mathbf{X}}_{k}} \right)$ and $\lambda \left( {{{\mathbf{\hat{X}}}}_{k}} \right)$, corresponding to the matrices involved in the SIAI and the SIAI\_F, respectively, during the 7 iterations are both listed in Table 3. The upper bound of the real bandwidth of ${{\mathbf{X}}_{k}}$ evaluated in terms of Theorem 5, $\frac{{{3}^{k}}+1}{2}\lambda \left( \mathbf{A} \right)$, is also listed in the last line of Table \ref{Tab2}. It can be seen from Table \ref{Tab2}, $\lambda \left( {{\mathbf{X}}_{k}} \right)$ and $\frac{{{3}^{k}}+1}{2}\lambda \left( \mathbf{A} \right)$ are the same in the first five iterative steps, and $\lambda \left( {{\mathbf{X}}_{k}} \right)$ is smaller than $\frac{{{3}^{k}}+1}{2}\lambda \left( \mathbf{A} \right)$ in the last two iterative steps. This is because many elements in ${{\mathbf{X}}_{6}}$ and ${{\mathbf{X}}_{7}}$ are so close to zero that the computer treated them as zeros. Additionally, $\lambda \left( {{{\mathbf{\hat{X}}}}_{k}} \right)$ is the same as $\lambda \left( {{\mathbf{X}}_{k}} \right)$ in the first two steps, and is far smaller than $\lambda \left( {{\mathbf{X}}_{k}} \right)$ in the last five steps, which means the filtering method improve significantly the sparsity of matrices during the iterative process. The numerical observation is consistent with the prediction of the theorem in Section \ref{S3}.

\begin{center}
\abovecaptionskip 0pt \belowcaptionskip 1pt
\renewcommand{\arraystretch}{1.2}
{\scriptsize \tabcaption{Comparisons of $\lambda \left( {{{\mathbf{\hat{X}}}}_{k}} \right)$, $\lambda \left( {{\mathbf{X}}_{k}} \right)$ and $\frac{{{3}^{k}}+1}{2}\lambda \left( \mathbf{A} \right)$}\label{Tab2}
\noindent\begin{tabular*}{\textwidth}{@{\extracolsep{\fill}}@{~~}cccccccc}
\toprule
{$k$} & {1} & {2} & {3} & {4} & {5} & {6} & {7}\\
\midrule
$\lambda \left( {{{\mathbf{\hat{X}}}}_{k}} \right)$ & 2 & 4 & 8 & 10 & 10 & 10 & 8\\
$\lambda \left( {{\mathbf{X}}_{k}} \right)$ & 2	& 4 & 10 & 28 &	82 & 206 & 232\\
$\frac{{{3}^{k}}+1}{2}\lambda \left( \mathbf{A} \right)$ & 2	& 4	& 10 & 28 & 82	&244 & 730\\
\bottomrule
\end{tabular*}\small }\\[4mm]
\end{center}

The computational times and the errors corresponding to the SIAI\_F, the SIAI and the sqrtm function are listed in Table \ref{Tab3}. It can be seen that the SIAI\_F performs much better than the SIAI and the sqrtm function in both computational efficiency and accuracy, which reflects the potential of the SIAI\_F in the computation of the matrix square roots with sparse approximations.

\begin{center}
\abovecaptionskip 0pt \belowcaptionskip 1pt
\renewcommand{\arraystretch}{1.2}
{\scriptsize \tabcaption{Comparisons of SIAI\_F, SIAI and sqrtm}\label{Tab3}
\noindent\begin{tabular*}{\textwidth}{@{\extracolsep{\fill}}@{~~}cccccc}
\toprule
\multicolumn{2}{c}{SIAI\_F}& \multicolumn{2}{c}{SIAI} & \multicolumn{2}{c}{sqrtm}\\
\cline{1-6} $er$ & Time (s) & $er$ & Time (s) & $er$ & Time (s)\\
\midrule
7.62E-15&	0.12& 	1.86E-15&	75.44& 	9.94E-13&	14.93\\
\bottomrule
\end{tabular*}\small }\\[4mm]
\end{center}

\subsection{Performance of the SIAI\_F}\label{S5.3}
In this experiment, 45 sparse adjacency matrices from https://networkrepository.com are used to test the performance of the SIAI\_F. The dimensions of these matrices range from 10000 to 343791. As it cannot be guaranteed that the adjacency matrix has a principal square root, we actually use the matrix ${{\mathbf{A}}_{i}}=\mathbf{I}-{0.5{{\mathbf{B}}_{i}}}/{\rho \left( {{\mathbf{B}}_{i}} \right)}\;$, where ${{\mathbf{B}}_{i}}$ is the adjacency matrix. The SIAI\_F with the tolerance relative error ${{\varepsilon }_{\text{tol}}}={{10}^{-14}}$ is used for the computation of these matrices. To test the computational efficiency of the SIAI\_F, the sqrtm function is also used here. As these matrices are too large, the DB, the IN, the iterative scheme (\ref{E9}), and the SIAI are not suitable. The sqrtm is only suitable to the first 14 matrices in all the 45 matrices, and fails for the last 31 matrices due to the insufficient memory. 

The computational times and errors comparisons between the SIAI\_F and the sqrtm are displayed in Figure 3(a), and 3(b), respectively. It can be observed that the SIAI\_F performs better than the sqrtm in terms of both the computational accuracy and efficiency.

\begin{figure}
\mbox
{\subfigure[Computational error]{\includegraphics[width=70mm]{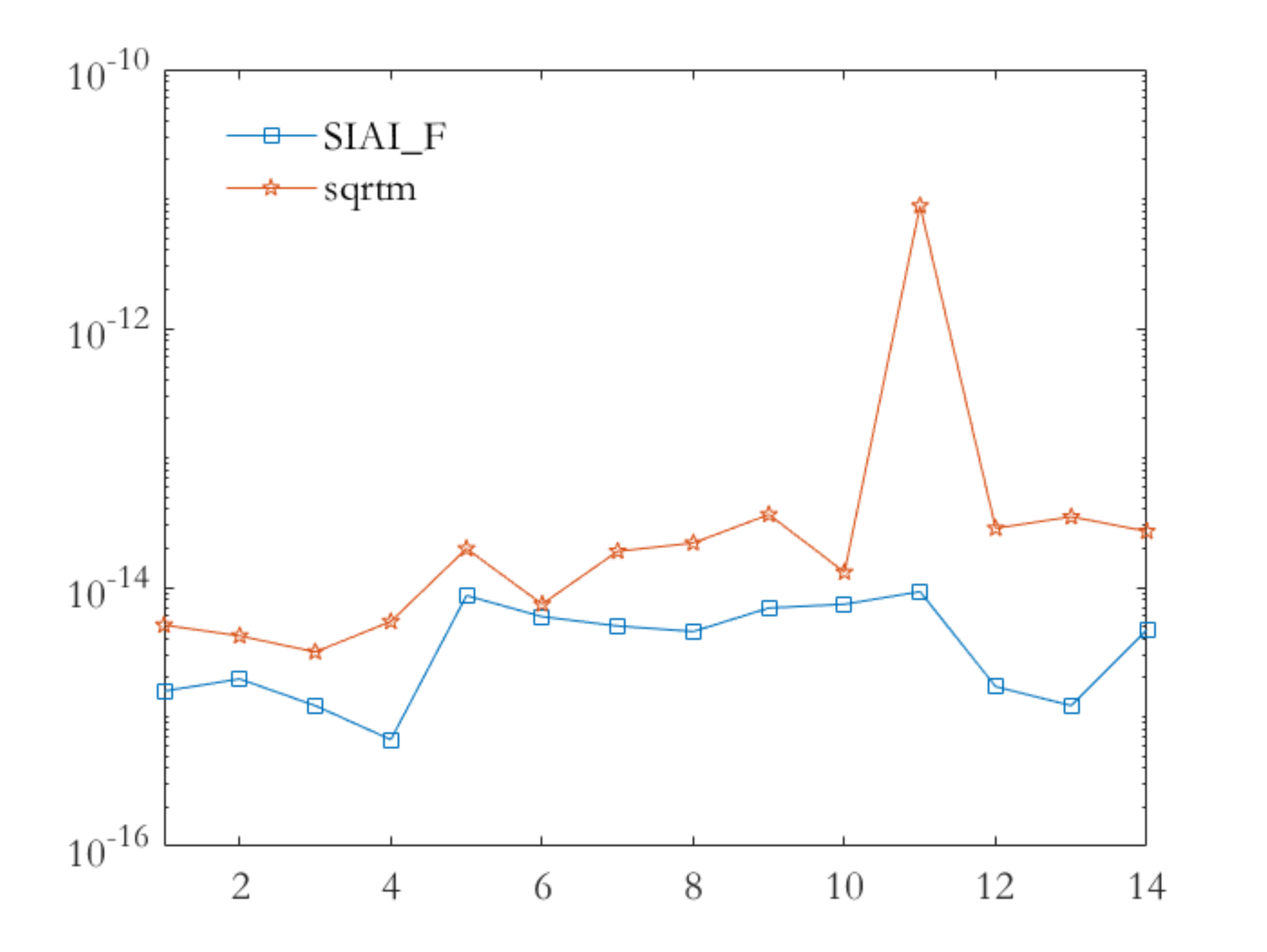}}
\quad 
\subfigure[Computational time]{\includegraphics[width=70mm]{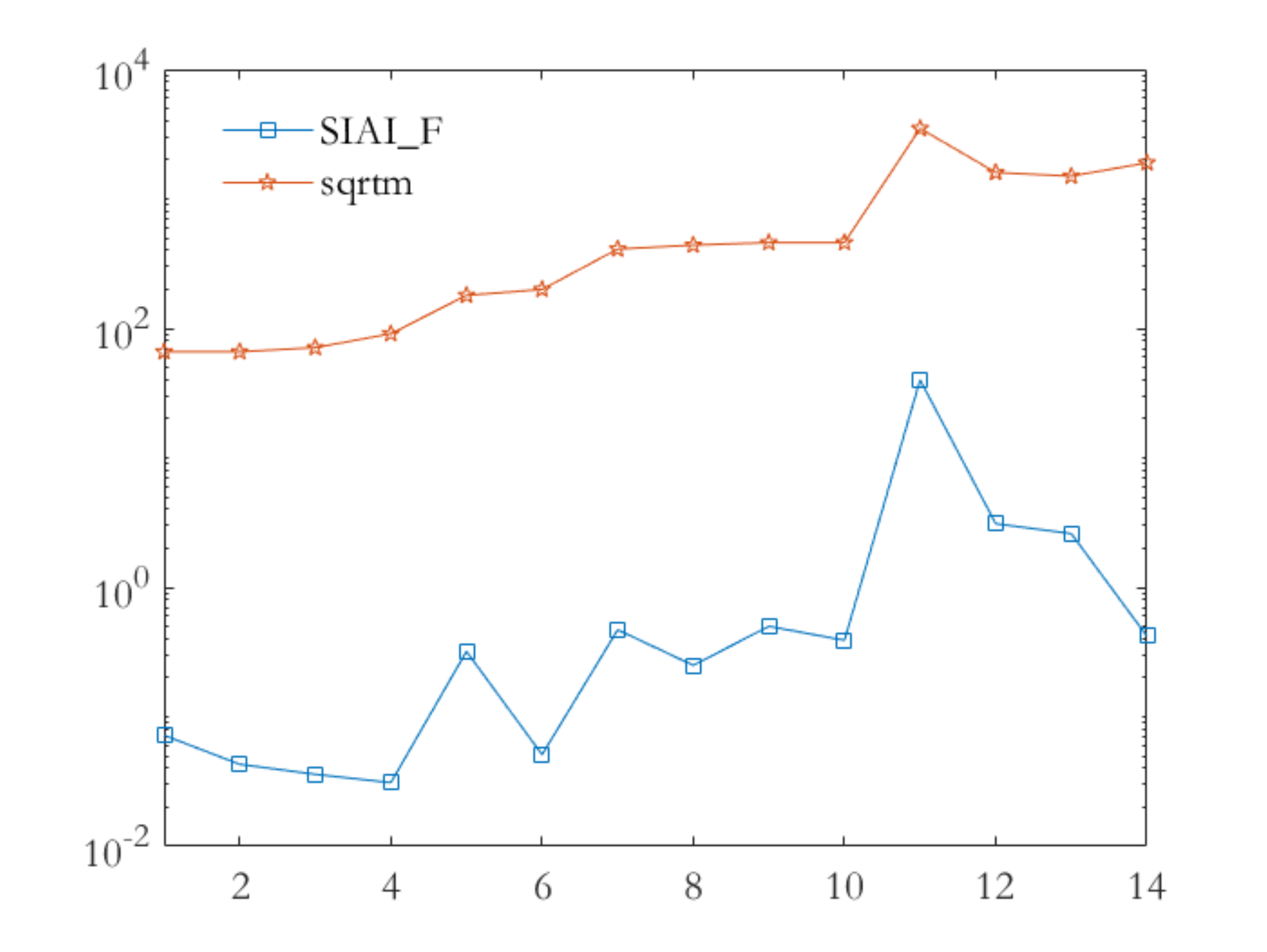}}} 
\figcaption{Comparison of the SIAI\_F with the sqrtm on computation accuracy, the transverse axises of these two figures are both the matrix number ${{\mathbf{A}}_{1}}\sim{{\mathbf{A}}_{14}}$, the longitudinal axis of (a) is the computational error, and the longitudinal axis of (b) is the computational time}\label{F3}
\end{figure}

For the last 31 matrices, the sqrtm fails due to the problem of insufficient memory. For these large and sparse matrices, the proposed SIAI\_F still work well. In Table 4, the dimension and sparsity of the 31 matrices are listed in columns 2 and 3. And the relative errors of the principal matrix square roots obtained by SIAI\_F and the time spent in the computation are listed in columns 4 and 5. The sparsity of ${{\mathbf{\hat{X}}}_{k}}$ is listed in the last column.

From Table \ref{Tab4} we can find that the SIAI\_F can compute the square roots of large-scale sparse matrices efficiently and accurately. Due to the using of adaptive filtering, the obtained principal square root matrix is also sparse, indicating that the SIAI\_F is very suitable for the computation when the $\varepsilon $-bandwidth of $\sqrt{\mathbf{A}}$ is much smaller than the real bandwidth.

\begin{table}[H]
\abovecaptionskip 0pt \belowcaptionskip 1pt
\renewcommand{\arraystretch}{1.2}
{\scriptsize
\tabcaption{Experiment results of ${{\mathbf{A}}_{15}}\sim{{\mathbf{A}}_{45}}$}\label{Tab4}
\noindent \begin{tabular*}{\textwidth}{@{\extracolsep{\fill}}@{~~}cccccc}
\toprule
No. & Dimension & Sparsity ($\mathbf{A}$) & $er$ & Time (s) & Sparsity (${{\mathbf{\hat{X}}}_{k}}$)\\
\midrule
        ${{\mathbf{A}}_{15}}$ & 38000 & 1.22E-04  & 8.02e-15 & 3.20E+00  & 2.30E-03  \\   
        ${{\mathbf{A}}_{16}}$ & 38120 & 8.35E-04  & 1.06e-14 & 7.90E+01  & 4.50E-03  \\  
        ${{\mathbf{A}}_{17}}$ & 38744 & 1.20E-03  & 7.43e-15 & 1.40E+01  & 3.70E-03  \\   
        ${{\mathbf{A}}_{18}}$ & 43471 & 1.09E-04  & 6.39e-15 & 2.70E+00  & 1.30E-03  \\   
        ${{\mathbf{A}}_{19}}$ & 43644 & 1.40E-04  & 2.89e-15 & 3.00E+00  & 1.80E-03  \\   
        ${{\mathbf{A}}_{20}}$ & 43747 & 1.09E-04  & 6.39e-15 & 2.60E+00  & 1.30E-03  \\   
        ${{\mathbf{A}}_{21}}$ & 56468 & 9.68E-05  & 6.11e-15 & 1.30E+02  & 8.50E-03  \\   
        ${{\mathbf{A}}_{22}}$ & 57735 & 9.98E-05  & 2.18e-15 & 5.50E-01  & 9.95E-05  \\   
        ${{\mathbf{A}}_{23}}$ & 68121 & 5.14E-04  & 6.17e-15 & 1.50E+00  & 1.89E-04  \\   
        ${{\mathbf{A}}_{24}}$ & 83995 & 7.19E-05  & 1.43e-15 & 1.50E+00  & 3.53E-04  \\   
        ${{\mathbf{A}}_{25}}$ & 91813 & 2.58E-05  & 7.92e-15 & 5.60E-01  & 7.87E-05  \\   
        ${{\mathbf{A}}_{26}}$ & 116670 & 2.60E-05  & 5.07e-15 & 1.50E+00  & 1.69E-04  \\   
        ${{\mathbf{A}}_{27}}$ & 122494 & 2.58E-05  & 5.33e-15 & 2.40E+00  & 2.66E-04  \\   
        ${{\mathbf{A}}_{28}}$ & 122750 & 2.30E-05  & 7.56e-15 & 1.70E+00  & 1.49E-04  \\   
        ${{\mathbf{A}}_{29}}$ & 126483 & 2.40E-05  & 4.72e-15 & 1.80E+00  & 1.73E-04  \\   
        ${{\mathbf{A}}_{30}}$ & 127998 & 2.37E-05  & 6.52e-15 & 1.80E+00  & 1.72E-04  \\   
        ${{\mathbf{A}}_{31}}$ & 131488 & 2.30E-05  & 7.19e-15 & 2.70E+00  & 2.67E-04  \\   
        ${{\mathbf{A}}_{32}}$ & 131490 & 1.98E-05  & 3.72e-15 & 2.10E+00  & 1.59E-04  \\   
        ${{\mathbf{A}}_{33}}$ & 135276 & 2.25E-05  & 6.83e-15 & 2.00E+00  & 1.61E-04  \\   
        ${{\mathbf{A}}_{34}}$ & 136239 & 2.23E-05  & 4.60e-15 & 1.80E+00  & 1.48E-04  \\   
        ${{\mathbf{A}}_{35}}$ & 140999 & 2.15E-05  & 5.16e-15 & 2.10E+00  & 1.51E-04  \\   
        ${{\mathbf{A}}_{36}}$ & 147772 & 2.06E-05  & 5.47e-15 & 2.30E+00  & 1.52E-04  \\   
        ${{\mathbf{A}}_{37}}$ & 152773 & 1.99E-05  & 6.24e-15 & 2.40E+00  & 1.44E-04  \\   
        ${{\mathbf{A}}_{38}}$ & 153563 & 1.98E-05  & 5.60e-15 & 2.30E+00  & 1.45E-04  \\   
        ${{\mathbf{A}}_{39}}$ & 158058 & 1.93E-05  & 5.77e-15 & 2.40E+00  & 1.42E-04  \\   
        ${{\mathbf{A}}_{40}}$ & 161070 & 7.82E-05  & 7.46e-15 & 1.20E+01  & 2.03E-04  \\   
        ${{\mathbf{A}}_{41}}$ & 162650 & 1.45E-05  & 7.94e-15 & 1.50E+00  & 7.52E-05  \\   
        ${{\mathbf{A}}_{42}}$ & 172175 & 1.77E-05  & 5.69e-15 & 2.60E+00  & 1.32E-04  \\   
        ${{\mathbf{A}}_{43}}$ & 203954 & 2.33E-05  & 2.40e-15 & 2.50E+00  & 8.55E-05  \\   
        ${{\mathbf{A}}_{44}}$ & 245874 & 5.51E-06  & 6.49e-15 & 6.00E-01  & 1.07E-05  \\   
        ${{\mathbf{A}}_{45}}$ & 343791 & 7.29E-06  & 9.04e-15 & 1.40E+00  & 1.58E-05  \\   
\bottomrule
\end{tabular*}\small }\\[4mm]
\end{table}

\section{Conclusion}\label{S6}
In this study, a new stable iterative method avoiding fully matrix inversions (denoted by SIAI) is proposed based on the Newton iteration. The corresponding proof on the numerical stability of the proposed method is given. The selection of the initial matrix for the iteration is also recommended, which broadens the applicable matrix range of the original Newton iteration. The sparsity of the principal matrix square root is also discussed by using the real bandwidth and $\varepsilon \text{-}$bandwidth. The analysis shows that there exists a sparser approximate matrix with the similar accuracy for each matrix involved in every iterative step of the SIAI. The filtering technique with an adaptive filtering threshold based on the error analysis is proposed to obtain this sparser approximate matrix. Combining the filtering technique with the SIAI yields the SIAI\_F. Numerical examples shows the contributed SIAI is more efficient than the Newton iteration, the DB method and the IN method. The proposed SIAI\_F can compute the matrix whose principal square root is nearly sparse efficiently and precisely. 

According to the analysis in Section 2, the efficiency of the proposed algorithm will be dependent on the $\varepsilon \text{-}$bandwidth of $\sqrt{\mathbf{A}}$. If the $\varepsilon \text{-}$bandwidth of $\sqrt{\mathbf{A}}$ is not far smaller than the dimension of the input matrix, the filtering technique can just do a very limited improvement in computational efficiency. However, we believe it is a valuable attempt for the computation of the square roots of sparse matrices. In the future, we will extend the proposed algorithm to the computation of the matrix $p$-th root.

\end{CJK*}
\end{document}